\RequirePackage{booktabs}
\documentclass[pdflatex,sn-mathphys,Numbered]{sn-jnl}

\usepackage[english]{babel}
\usepackage[T1]{fontenc}
\usepackage[utf8]{inputenc}

\usepackage{graphicx}%
\usepackage{multirow}%
\usepackage{amsmath,amssymb,amsfonts}%
\usepackage{amsthm}%
\usepackage{mathrsfs}%
\usepackage[title]{appendix}%
\usepackage[table]{xcolor}%
\usepackage{textcomp}%
\usepackage{manyfoot}%
\usepackage{booktabs}%
\usepackage{algorithm}%
\usepackage{algorithmicx}%
\usepackage{algpseudocode}%
\usepackage{listings}%

\usepackage{array}
\usepackage{tabularx}
\usepackage{csquotes}
\usepackage{alltt}
\usepackage{colortbl}
\usepackage[caption=false]{subfig}

\usepackage{placeins}

\newcommand{\CP}{C\&P}

\newcommand{\NN}{\mathbb{N}}

\newtheorem{example}{Example}[section]

\newtheorem{definition}{Definition}[section]

\begin{document}

\title{Constraint programming methods in three-dimensional container packing}

\author{\fnm{Szymon} \sur{Wróbel}}
\email{szymon.wrobel@pwr.edu.pl}

\affil{
    \orgdiv{Department of Fundamentals of Computer Science} \\
    \orgname{Wrocław University of Science and Technology}
}

\abstract{
Cutting and packing problems are present in many, at first glance unconnected, areas, therefore it's beneficial to have a good understanding of their underlying structure, to select proper techniques for finding solutions.

Cutting and packing problems are a class of combinatorial problems in which there are specified two classes of objects: big and small items and the task is to place the small items within big items.

Even in the 1-dimensional case, bin-packing is strongly NP-hard (see \cite{Garey1978}), which suggests, that exact solutions may not be found in a reasonable time for bigger instances. In the literature, there are presented many various approaches to packing problems, e.g. mixed-integer programming, approximation algorithms, heuristic solutions, and local search algorithms, including metaheuristic approaches like Tabu Search or Simulated Annealing.

The main goal of this work is to review existing solutions, survey the variants arising from the industry applications, present a solution based on constraint programming and compare its performance with the results in the literature. Optimization with constraint programming is a method searching for the global optima, hence it may require a higher workload compared to the heuristic and local search approaches, which may finish in a local optimum.

The performance of the presented model will be measured on test data used in the literature, which were used in many articles presenting a variety of approaches to three-dimensional container packing, which will allow us to compare the efficiency of the constraint programming model with other methods used in the operational research.
}


\maketitle

\section{Introduction}
\label{ch:introduction}

Cutting and packing problems (\CP) are a family of combinatorial problems, that appear in many, at first glance unconnected, areas -- from patterns in cutting elements from bigger pieces, through minimization of the number of used containers needed for delivery, to allocation of resources and task scheduling.

Cutting and packing problems appear in many different variations, but all of them exhibit certain common properties. The first one is a distinction between two groups of items: small and big ones\footnotemark -- big items are containers, collective packages etc., while small items are subject to packaging. The next common property is satisfying the \textit{geometric constraint}, meaning the items may not overlap and small items are fully inside the containers. The last commonality is, that the expected result is the packing pattern, that optimizes some objective function -- most often maximization of volume utilization, minimization of height of packing, or minimization of number of containers used.

\footnotetext{In some works they are called objects and items.}

Cutting and packing problems were identified and analyzed a long time ago, but a systematic approach to solving them started appearing in the 1960s (see P.~Sweeney, E.~Paternoster \cite{Sweeney1992}). Even though most of the scientific works are concerned with one- and two-dimensional packing, three-dimensional packing has many applications in industrial and business contexts. Having said that, results in the area of one- and two-dimensional variants have a significant impact on three- and higher-dimensional problems, like in the case of  \textsc{Harmonic} heuristic \cite{Epstein2005}.

\subsection{Problem statement}

\CP~family consists of many common problems. One of the simplest is \textsc{Bin packing}. In the decision version of the problem (see \cite[problem SR1]{Garey1990}), the question is, if using $K$ containers of size $B$ it is possible to pack items $U = {u_1, \ldots, u_n}$ of sizes $w_1, \ldots, w_n \leq B$. Another way of looking at the problem is splitting the set $U$ into $K$ subsets, such that the sum of the sizes of the elements doesn't exceed $B$.

The optimization version of the problem is presented as follows (cf. \cite[rozdział~8]{Martello1990}):

\begin{align*}
    \min        \quad   & \sum_{k=1}^n y_k                              &               &   \\
    \text{s.t.} \quad   & \sum_{i=1}^n x_{k,i} \cdot w_i \leq B y_{k}   & k=1 \ldots n  &   \\
                        & \sum_{k=1}^M x_{k,i} = 1                      & i=1 \ldots n  &   \\
                        & x_{k,i} \in \{0,1\}                           & k=1 \ldots M, i=1 \ldots n & \\
                        & y_{k} \in \{0,1\}                             & k=1 \ldots M  &   \\
\end{align*}

\noindent where $y_k$ are decision variables describing usage of $k$-th container, and $x_{k_i}$ are binary variables describing placing $i$-th item in the $k$-th container.

% \begin{theorem}
% \textsc{Bin Packing} jest problemem NP-zupełnym
% \end{theorem}
% \begin{proof}
% Oczywiście, \textsc{Bin Packing} należy do klasy NP -- wybierając niedeterministycznie podzbiory możemy sprawdzić, czy suma wielkości elementów w~każdym z~podzbiorów nie przekracza pojemności.

% Aby pokazać NP-trudność problemu \textsc{Bin Packing}, możemy dokonać redukcji z~problemu \textsc{Partition}\footnotemark. Istotnie, niech $c_1, \ldots, c_n$ będą danymi dla problemu \textsc{Partition}. Oznaczmy przez $M$~sumę podanych liczb $\sum_{i=1}^n c_i$. Jeśli $M$~jest liczbą nieparzystą, to od razu wiadomo, że nie jest możliwy podział i~zwracamy odpowiedź \textsc{NIE}. W~przeciwnym wypadku wyznaczamy liczbę $K$, taką, że $M = 2K$.

% Zauważmy, że w~tym momencie pytanie o~podział liczb $c_1, \ldots, c_n$ jest równoważny zapytaniu, czy mając 2 kontenery o~pojemności $K$ możemy spakować przedmioty o~wielkościach $c_1, \ldots, c_n$. Wobec tego, jeśli uzyskamy odpowiedź dla tego egzemplarza \textsc{Bin Packing}, automatycznie uzyskujemy odpowiedź dla wejściowego egzemplarza \textsc{Partition}, ale \textsc{Partition} jest problemem NP-trudnym oraz w~redukcji użyliśmy jedynie pamięci logarytmicznej (liczby $M$ i~$K$), a~więc \textsc{Bin Packing} również jest NP-trudny.

% \footnotetext{\textit{pol.} \textsc{Podzieł}, problem NP-zupełny należący do listy problemów omawianej w~pracy R. Karpa \cite{Karp1972} -- Niech $c_1, \ldots, c_n \in \NN$. Czy istniej podzbiór $i~\subseteq \{1,\ldots,n\}$ t. że $\sum_{i~\in I} c_i = \sum_{i~\notin I} c_i$.}

% \end{proof}

Due to NP-completeness of \textsc{Bin Packing}, which was presented i.a. in \cite{Garey1978}, there are no polynomial algorithms solving that problem, unless $P = NP$. Additionally, in the same work, a stronger statement was proven: \textsc{Bin Packing} is a strongly NP-complete problem, which means there are no pseudopolynomial algorithms (runtime polynomial in size of the input values, not necessarily length of the input), contrary to \textsc{Knapsack} problem, where there exists a dynamic programming algorithm working in $\mathcal{O}(nB)$ time.

Naturally, NP-completeness of higher dimensional variants of \textsc{Bin Packing} is a trivial consequence of NP-completeness of one-dimensional version. A reduction for \textsc{$d$-D Bin Packing} is created by fixing dimensions $2, \ldots, d$ for containers and items, that is for input data for \textsc{Bin Packing}:
\begin{align*}
    \text{items}~c_1, \ldots, c_n \in \NN, \\
    K \text{containers of volume}~B
\end{align*}

the input data for \textsc{$d$-D Bin Packing} would be:

\begin{align*}    
    \text{items}~(c_1, 1, \ldots, 1), \ldots (c_n, 1, \ldots, 1) \in \NN^d, \\
    K \text{containers of dimensions}~(B, 1, \ldots, 1)
\end{align*}

\subsection{\CP~Problem Typology}
\label{sec:typology}

Cutting and packing problems appear in many research areas, therefore it's important to create a way to systematically describe the key properties of the analyzed problem. The first broader known such system is a typology introduced by H.~Dyckhoff in 1990~\cite{Dyckhoff1990}.

\subsubsection{Dyckhoff Typology}
\label{sec:dyckhoff}
Dyckhoff's typology focuses on four criteria, that make up the problem description. First of them is dimensionality -- the number of dimensions needed to describe a packing pattern. The second criterion is kind of assignment, which means if we are packing all small items and selecting containers (\textit{ger.} Verladeproblem), or if we are using all containers and picking small items (\textit{ger.} Beladeproblem)\footnotemark. The third aspect is the characterization of the containers: single~(O), many identical shapes~(I), and many different shapes~(D). The last criterion is the description of the small item assortment: few items of different shapes~(F), many items of many different shapes~(M), many items of a few different shapes~(R), many items of congruent shapes~(C). The summary of the criteria is presented in table \ref{tab:dyckhoff}.

\footnotetext{To the knowledge of the author there is no precise English term that can be used to differentiate between the two terms, as both would translate to ``loading problem''. The former is used in logistics contexts with regard to loading cargo onto transport vehicles (transshipment), while the latter is used for packing and efficient space usage.}

It's quite easy to notice, that within Dyckhoff's typology, there are 96 distinct classes of \CP~problems, but as it was pointed out in \cite{Washer2007}, sometimes it's not apparent to which class a problem would belong, for example, strip packing (problem of packing rectangles into a rectangle of least height) was placed in \texttt{2/V/O/M} class, while Dyckhoff classified it as \texttt{2/V/D/M}, which was supported by a remark, that the set of big items is a set of rectangles of all possible heights (cf. \cite{Dyckhoff1990}).

\begin{table}[!ht]
    \centering
    \begin{tabular}{p{10cm}}
    \begin{enumerate}
        \item Dimensionality
        \begin{description}
            \item[(1)] one-dimensional,
            \item[(2)] two-dimensional,
            \item[(3)] three-dimensional,
            \item[(N)] N-dimensional ($N > 3$).  
        \end{description}

        \item Kind of assignment
        \begin{description}
            \item[(B)] all containers, selection of small items,
            \item[(V)] selection of containers, all small items.
        \end{description}

        \item Assortment of large objects
        \begin{description}
            \item[(O)] single container,
            \item[(I)] many identical containers,
            \item[(D)] many different containers.
        \end{description}

        \item Assortment of small items
        \begin{description}
            \item[(F)] few items of different shapes,
            \item[(M)] many items of different shapes,
            \item[(R)] many items of a relatively few shapes,
            \item[(C)] items of congruent shapes.
        \end{description}
    \end{enumerate}
    \\
    \end{tabular}
    \caption{Dyckhoff's Typology (based on \cite{Dyckhoff1990})}
    \label{tab:dyckhoff}
\end{table}

\subsubsection{W\"ascher Typology}

The Dyckhoff typology, even if quite extensive, isn't perfect. After noticing some issues with it, W\"asher \cite{Washer2007} presented their own classification of \CP~problems. The Main criteria introduced by W\"asher were: dimensionality, kind of assignment, assortment of big items, assortment of small items, and shape of small items.

The kind of assignment is a classification criterion analogous to the \textit{Beladen/Verladen} dichotomy, but to avoid issues with comprehending the difference for people not speaking German (cf. footnote in section \ref{sec:dyckhoff}), a new division was introduced: minimization of the cost and maximization of value. The next property describing \CP~problems is the assortment of small items: identical, weakly heterogeneous, and strongly heterogeneous. Weakly heterogeneous assortment consists of a small (relatively to number of items) number of classes, while strongly heterogeneous assortment consists of items with heavily diversified shapes and sizes -- there are not many similar items.

The kind of assignment together with an assortment of small items constitute the basic problem types. Taking the assortment of big items leads to intermediate problem types. Categories appearing in the analysis of that criterion are the number of the big items and their character. The first group consists of problems with a single big item: it may have fixed dimensions, or at least one of the dimensions may be open (potentially unbounded, but may be part of the objective function). The second group contains problems with multiple big items, whose assortment may be, similarly to small items, identical, weakly- and strongly heterogeneous.

Dimensionality is treated in a similar way as in Dyckhoff's typology, i.e. number of dimensions of concern for the problem. The last criterion is the shape of small items -- regular shapes (rectangles, boxes, cylinders, ...), or irregular. Including those two properties, results in refined/combined problem types.

\paragraph{Characterization of basic problem types}

Considering the criteria of the kind of assignment and assortment of small items leads to basic problem types: \textsc{Identical Item Packing Problem}, \textsc{Placement Problem}, \textsc{Knapsack Problem}, \textsc{Open Dimension Problem}, \textsc{Cutting Stock Problem}, and \textsc{Bin Packing Problem}.

The first three problems are maximization problems, differing in the small item assortment: identical items -- \textsc{Identical Item Packing Problem} (IIPP), weakly heterogeneous -- \textsc{Placement Problem} (PP), strongly heterogeneous -- \textsc{Knapsack Problem} (KP). A common properties of maximization problems are limited number of big items, and very often excess of small items, which leads to a need for a selection of items that will be packed.

The remaining problems are minimization problems, where the objective is to pack all small items while minimizing the number of required big items. In the case of \textsc{Open Dimension Problem} (ODP) one or more of the dimensions of the big item are variable, and the objective is minimization of its size. In the case of \textsc{Cutting Stock Problem} (CSP) and \textsc{Cutting Stock Problem} (CSP), the dimensions of big items are fixed.

\subsection{Problem variants and extensions}
\label{sec:rozszerzenia}

In the W\"ascher typology, a number of implicit assumptions were introduced, for example, single objective function, orthogonal packing pattern, regular shape of the items, and upfront knowledge of the payload (offline problems). It's readily apparent, that the problems arising from business applications do not always allow for such assumptions, for example, packing of L-shaped items (cf. Wong, Yung \cite{Wong2010}), irregular shape of the container (air freight containers, cf. Paquay et al. \cite{Paquay2016}), or the assortment of the items is not known upfront -- payload is arriving dynamically (cf. Hemminki et al. \cite{Hemminki1998}).

Outside of the implicit assumptions, the primary assumption in the problems was a geometrical constraint, i.e. the small items don't overlap and are fully contained within the big item(s). One of the most often appearing extensions of the problem is allowing item rotation -- in many cases, there are no contraindications to rotating the items, but in some cases, the vertical orientation must be preserved (e.g. fragile payload, ``This side up!'') considered i.a. in L.~Epstein, R.~van Stee \cite{Epstein2006}.

As might be expected, in industrial applications, additional physical constraints may arise: stability (vertical, horizontal, longitudinal), and load distribution constraints (cf. guidelines on load securing \cite{IRU2014}). Very common extensions in literature are also multi-drop constraint (many delivery points) and load bearing constraint. The constraints arising from practical applications are the subject of many works and surveys, i.a. Liu et al. \cite{Liu2011}, Bischoff, Ratcliff \cite{Bischoff1995}, Junqueira et al. \cite{Junqueira2012a,Junqueira2012, Junqueira2013}.

\section{Selected Topics of Constraint Programming}
\label{ch:constraint}

Constraint programming is a declarative programming paradigm, in which the problem is modelled using the variables, with their domains (set of possible values) and constraints describing relations between the variables. Many of the common problems contain some constraints, e.g. lecture B must start after the end of lecture A. A typical example of a problem described with constraints is the $n$-queen problem, which poses a task of placing $n$ queens on the $n \times n$ board in such a way, that no two pieces threaten itself.

\subsection{Constraint satisfaction problem}

The Constraint Satisfaction Problem (CSP) is the underlying problem of constraint programming. Formally, it can be defined as follows (based on \cite{Rossi2006}):

\vspace{1em}
\begin{definition}
\textbf{Constraint Satisfaction Problem} is a triple $\mathcal{P} = (\mathcal{X}, \mathcal{D}, \mathcal{C})$, where $\mathcal{X} = (x_1, x_2, \ldots, x_n)$ is a set of variables, $\mathcal{D} = (D_1, D_2, \ldots, D_n)$ a set of domains of the variables -- $\textbf{dom}\;x_i = D_i$, and $\mathcal{C} = (C_1, C_2, \ldots, C_k)$ a set of constraints. Constraint is a pair $C_i = (S_i, R_i)$, such that, $S_i \subseteq \mathcal{X}$ is a subset of the variables, and $R_i \subseteq \prod_{x \in S_i} \textbf{dom}\;x$ is a relation defined over the domains of the subset of the variables $S_i$.

A solution of the problem $\mathcal{P}$ is an assignment $A = (a_1, \ldots, a_n)$ such that $a_i$ is the value of variable $x_i$, $a_i \in D_i$, and for each constraint $C_j \in \mathcal{C}$, a projection of $A$ on $S_j$ is in relation $R_j$
\end{definition}
\vspace{1em}

With the definition above, in examples \ref{ex:n-queens} and \ref{ex:send-more-money} presented are classical constraint satisfaction problems, the previously mentioned $n$-queen problem, and the cryptarithmetic puzzle \texttt{SEND + MORE = MONEY}.

\vspace{1em}
\begin{example}[$n$-queens problem]
\label{ex:n-queens}

Let $\mathcal{X} = (x_1, \ldots, x_n)$ be a set of variables describing the row in which a queen is placed in the $i$-th column. The domains of the variables $\mathcal{D} = (D_1, \ldots, D_n)$ are equal, $D_i = \{1,\ldots, n\}$. For each pair of columns $i, j$, two constraints are present: 

\begin{itemize}
    \item $x_i \neq x_j$ -- queens don't threaten each other along the row,
    \item $\left| x_i - x_j \right| \neq \left|j - i\right| $ -- queens don't threaten each other along the diagonals.
\end{itemize}

The last remaining constraint arising from the description -- non-threatening in the column -- is trivially satisfied by the careful selection of the meaning of variables, based on a simple observation, that in each column there must be exactly one queen -- otherwise they would threaten each other.

\end{example}

\vspace{1em}

\begin{example}[Cryptarithmetic puzzle]~\\
\label{ex:send-more-money}

The objective of the puzzle is to substitute letters with digits in such a way, that the equation is satisfied, the same digit corresponds to the same letters, and each letter corresponds to a different digit. 

\begin{center}
\begin{tabular}{*{5}{c}}
    & S & E & N & D \\
  + & M & O & R & E \\
  \hline
  M & O & N & E & Y \\
\end{tabular}
\end{center}

The set of variables consists of variables corresponding to the letters:

\[ \mathcal{X} = (x_S, x_E, x_N, x_D, x_M, x_O, x_R, x_Y)\]

Values of the variables are digits, therefore

\[ x_S, x_E, x_N, x_D, x_M, x_O, x_R, x_Y \in \{0, \ldots, 9\} \]

Letters $S$ and $M$ appear as a leading digit, therefore they can't be equal to zero:

\[ x_S \neq 0, \quad x_M \neq 0 \]

Each letter corresponds to a different digit, therefore:

\[ i \neq j \implies x_i \neq x_j \]

The last constraint describes satisfying the equation
\[ 
    x_M \cdot 10^4 + 
    (x_O - x_S - x_M) \cdot 10^3 + 
    (x_N - x_E - x_O) \cdot 10^2 + 
    (x_E - x_N - x_R) \cdot 10 + 
    (x_Y - x_D - x_E) 
    = 0 
\]

\end{example}
\vspace{1em}

\subsection{Local consistency}

The first trivial approach to solving CSP is conducting a full search of all of the possible assignments, for example, using backtracking. Unfortunately, very often occurs an effect called ``thrashing'', which is a repeated assignment of values leading to breaking some constraints, without changing the assignments impacting the satisfiability in the subtree (cf. \cite[ch. 2.2.3]{Rossi2006}). This is a consequence of the fact, that in the basic version, the backtracking algorithm doesn't ``learn'', that is, doesn't use the information from prior assignments leading to backtracking. One of the approaches to improving the learning process, and reducing the number of wasteful backtracking are methods of achieving local consistency, also called constraint propagation algorithms.

\subsubsection{Node consistency}

The simplest notion of local consistency is node consistency. An unary constraint is said to be node-consistent if all elements of the domain satisfy that constraint; the problem is node-consistent if all unary constraints are node-consistent.

An approach to node consistency can be presented in the example \ref{ex:send-more-money}: constraint $x_S \neq 0$ is not node consistent, because the value $0$ from the domain of $x_S$ ($D_S = \{0, \ldots, 9\}$) doesn't satisfy that constraint. Naturally, node consistency can be achieved by removing the value $0$ from the domain.

\subsubsection{Arc consistency}

The next concept of local consistency is arc consistency. A binary constraint $C = ((x_i, x_j), R_{i,j})$ is said to be arc-consistent if for each value $x_i$ there exists a value $x_j$, such that $(x_i, x_j) \in R_{i,j}$, and symmetrically, for each $x_j$ there exists a value of $x_i$ that satisfying the constraint. Formally:

\[ 
    (\forall_{x_i \in D_i} \exists_{x_j \in D_j}, (x_i, x_j) \in R_{i,j}) 
    \land 
    (\forall_{x_j \in D_j} \exists_{x_i \in D_i}, (x_i, x_j) \in R_{i,j}) 
\]

The simplest algorithm for enforcing the arc consistency is removing those values from domains, that don't have support, i.e. there are no values in the other domain that would satisfy the constraint with it. This algorithm is known as \textsc{AC-1}, but there are more sophisticated algorithms for enforcing arc consistency (cf. \cite{Mackworth1977}).

An important observation is, that the arc consistency doesn't mean there exists a solution, which was presented in the example \ref{ex:arc-consistency-failure}.

\vspace{1em}
\begin{example}[Arc consistency and satisfiability]
\label{ex:arc-consistency-failure}
Let $x_1, x_2, x_3 \in \{0,1\}$ be boolean variables satisfying constraints $x_1 \neq x_2$, $x_2 \neq x_3$, $x_1 \neq x_3$. This problem is obviously arc-consistent but has no solution.
\end{example}
\vspace{1em}

\subsubsection{Other local consistency concepts}

Even though node and arc consistency often allow for a significant reduction of the search space, they are fundamental notions of local consistency. In the literature presented were various other concepts of consistency, e.g. path consistency, or hyper-arc consistency -- an extension of arc consistency to $k$-ary constraints (cf. \cite{Apt2003}, \cite{Rossi2006}).

\subsection{Global Constraints}
\label{sec:wiezy-globalne}

Global constraints are special constraints that describe a specific relationship between variables, for example, constraint \textsc{AllDifferent}$(x_1,\ldots,x_n)$ means, that the variables $x_1, \ldots, x_n$ are pairwise different. Of course, it is possible to describe it using simpler basic constraints, but a common programming practice is \textsc{DRY} -- Don't Repeat Yourself, which advises against repeating the same patterns. Global constraints can also improve the readability and semantic meaning of the model, as well as reduce the number of points, where a designer can make an error.

Besides pragmatic reasons, global constraints can also allow the solver to use specialized propagation and search algorithms, thanks to the known structure of the constraint, instead of the default generic methods. (cf. \cite{Apt2003}). Global constraints can also be treated as idiomatic expressions -- sources such as \cite{GCC2010} are cataloguing and describing the most common global constraints, and are presenting their properties and limitations.

\section{Literature review}
\label{ch:literature}

Cutting and packing problems appear in many branches of industry, where they are used in various manners, mostly for optimizing logistic tasks, as well as scheduling, waste minimization or planning. In industrial applications, reaching global optimum is usually not as important, as getting a satisfying result quickly, which leads to a high popularity of heuristic, approximation or metaheuristic approaches.

A summary of the surveyed literature, regarding cutting and packing problems is presented in the table \ref{tab:literatura}, with a description of the tasks considered in the works and presented approaches. The first three columns describe the type of problem: ODP (Open Dimension Problem), BPP (Bin Packing Problem), SCP (Single Container Packing). In the next column presented are variants and extensions considered in the work: OL (online problem), ROT (rotation of objects), STAB (load stability criterion), MD (Multi-Drop Constraint). Last group of columns describes the approach selected by the authors: MIP (Mixed-Integer Programming), HEU (heuristic), MH (metaheuristic), APR (approximation algorithm).

\begin{table}[!ht]
\centering
\rowcolors{5}{white}{black!10}
\begin{tabular}{  >{\raggedright\arraybackslash}m{4cm} cccc cccc }
    \toprule
    &
        \multicolumn{3}{c}{Type} & 
        \multirow{2}{*}{Variant} &
        \multicolumn{4}{c}{Approach} \\
    \cmidrule{2-4}
    \cmidrule{6-9}
    & 
        ODP & BPP & SLOPP & &
        MIP & HEU & MH & APR \\
    \hline
    Li, Cheng \cite{Li1990} (1990) & 
        + & & & &
        & + & & + \\
    \hline
    Ngoi et al. \cite{Ngoi1994} (1994) &
        & & + & ROT &
        & & & \\
    \hline
    Bischoff et al. \cite{Bischoff1995a} (1995) &
        & & + & Survey &
        & + & & \\
    \hline
    Bischoff, Ratcliff \cite{Bischoff1995} (1995) &
        & & + & Survey &
        & + & & \\
    \hline
    Miyazawa et al.  \cite{Miyazawa1997} (1997) &
        + & & & &
        & + & & + \\
    \hline
    Bortfeldt, Gehring \cite{Bortfeldt1998, Bortfeldt1998a} (1998) &
        & & + & ROT &
        & & TS & \\
    \hline
    Hemminki et al. \cite{Hemminki1998} (1998) &
        & ? & & OL &
        & ? & & \\
    \hline
    Abdou, Elmasry \cite{Abdou1999} (1999) &
        & & + & OL &
        & + & & \\
    \hline
    Miyazawa et al. \cite{Miyazawa1999} (1999) &
        + & & & ROT &
        & + & & + \\
    \hline
    Eley \cite{Eley2002} (2002) &
        & & + & ROT &
        & + & & \\
    \hline
    Faroe et al. \cite{Faroe2003} (2003) &
        & + & & &
        & & GLS & \\
    \hline
    Jin et al. \cite{Jin2004} (2004) &
        & & + & ROT, STAB &
        & & SA & \\
    \hline
    Epstein, van Stee \cite{Epstein2005} (2005) &
        & + &&  OL &
        & + & & \\
    \hline
    Tsai, Li \cite{Tsai2006} (2006) &
        + & & & &
        + & & & \\
    \hline
    Epstein, van Stee \cite{Epstein2006} (2006) &
        + & + & & ROT &
        & + & & + \\
    \hline
    W\"asher et al. \cite{Washer2007} (2007) &
        \multicolumn{8}{ c }{C\&P problem typology } \\
    \hline
    Huang, He \cite{Huang2009} (2009) &
        & & + & &
        & + & & \\
    \hline
    Epstein, Levy \cite{Epstein2010} (2010) &
        & + & & ROT, OL\footnotemark[1] &
        & + & & + \\
    \hline
    He, Huang \cite{He2010} (2010) &
        & & + & &
        & + & & \\
    \hline
    Jiang, Cao \cite{Jiang2012} (2012) &
        & + & & ROT, STAB &
        & & SA & \\
    \hline
    Junqueira et al. \cite{Junqueira2012a} (2012) &
        & & + & MD &
        + & & & \\
    \hline
    Junqueira et al. \cite{Junqueira2012} (2012) &
        & & + & STAB, LB &
        + & & & \\
    \hline
    Bortfeldt, W\"ascher \cite{Bortfeldt2013} (2013) &
        \multicolumn{8}{ c }{Survey of problems and constraints} \\
    \hline
    Junqueira et al. \cite{Junqueira2013} (2013) &
        & & + & Survey &
        + & & & \\
    \hline
    Tsai et al. \cite{Tsai2015} (2015) &
        + & & & &
        + & & & \\
    \hline
    Toffolo et al. \cite{Toffolo2017} (2017) &
        & + & & STAB &
        & + & & \\
    \hline
    Huang, Hwang \cite{Huang2018} (2018) &
        + & & & &
        + & & & \\
    \hline
    Gzara et al. \cite{Gzara2020} (2020) &
        & + & & LB &
        +$^*$ & & & \\
    \hline
    \hiderowcolors
    \multicolumn{9}{ |l| }{}\\
    \multicolumn{9}{ |l| }{\textbf{Notations:}}\\
    \multicolumn{9}{ |l| }{
        \begin{tabular}{ c l }
            OL & Online formulation \\
            ROT & Object rotation \\
            MD & Multi-Drop Constraint \\
            LB & Load Bearing Constraint \\
        \end{tabular}
        \begin{tabular}{ c l }
            STAB & payload stability \\
            TS & Tabu Search \\
            GLS & Guided Local Search \\
            SA & Simulated Annealing \\
        \end{tabular}
    }\\
    \multicolumn{9}{ |l| }{\textbf{Remarks:}}\\
    \multicolumn{9}{ |l| }{
        \begin{tabular}{ c p{12cm} }
            \footnotemark[1] & dynamic packing -- items arrive and depart during runtime, the objective is the minimization of the number of used containers (see. \cite{Coffman1983})\\
        \end{tabular}
    }\\
    \multicolumn{9}{ |l| }{}\\
    \hline
\end{tabular}

\caption{Survey of constraints and solutions in literature}
\label{tab:literatura}
\end{table}

\subsection{Solution approaches}

Due to widespread applications of cutting and packing problems, diverse requirements on expected accuracy and time constraints, a wide array of approaches for finding solutions was presented in the literature.

\subsubsection*{Mathematical programming}

Among the methods for achieving global optimum, the most commonly found approach is mixed integer programming (MIP). The solutions are based on linear constraints and linear objective function, but contrary to linear programming (LP), some variables may be integer values. The most common approach to solving mixed integer programs is the branch-and-bound method, presented in 1960 by A. Land and A. Doig \cite{Land1960}.

\subsubsection*{Approximation algorithms}

Many if not most of the cutting and packing problems are NP-hard optimization problems, therefore unless $P = NP$, there are no polynomial algorithms for them. On the other hand, for many of the NP-hard optimization problems exist approximation algorithms, which may not provide an exact solution, but provide a bound on the error. For example, an $\epsilon$-approximation algorithm for a minimization problem provides a guarantee, that a solution is not higher than $\epsilon \cdot \texttt{OPT}$, where \texttt{OPT} is the global minimum.

\subsubsection*{Heuristic approach}

In many cases, finding the approximation constant of the algorithm is not a trivial task, but in average case the quality of the solutions is acceptable. In such a case, the algorithm is called a heuristic (gr. \textit{gr.} $\varepsilon \upsilon \rho \iota \sigma \kappa \omega $, eurisko -- find, discover). Such algorithms are very often used as subprocedures for other approaches -- for example, to find an initial point for local search algorithms).

\subsubsection*{Local search and metaheuristic algorithms}

The next, very broad branch of the methods are algorithms based on local search. In this approach, the search space consists of all the feasible solutions and for each one a neighbourhood -- a set of ``close'' solutions -- can be established. For example, in the travelling salesperson problem, the search space is the set of $n$-element permutations (order of visiting the cities), and the neighbourhood consists of all possible swaps with other cities.

Local search is the simplest, greedy approach, where in each step, the best solution from the neighbourhood is selected. It's easy to notice, that this method may get stuck in the local optimum, where every solution in the neighbourhood is worse than the current one, but in the whole space there may exist a better solution.

One of the possible modifications of that problem are metaheuristic algorithms, which are heuristics governing the search process. Good examples of such methods are Tabu Search (TS) (cf. F. Glover \cite{Glover1989, Glover1990}) and Simulated Annealing (SA) (cf. S. Kirkpatrick et al. \cite{Kirkpatrick1983}). The first method is based on remembering a few last moves, and eliminating them temporarily from the environment -- they become a ``tabu''. The other method is based on a physical phenomenon called annealing, in which metal is heated to a very high temperature and then slowly cooled down. In Simulated Annealing the temperature impacts the probability of picking a solution worse than current one, which allows for a further exploration of the search space and escaping the local minima.

\subsubsection*{Constraint programming}

The last of the considered approaches, the one presented in this work, is application of constraint programming techniques to the \CP~problems. Constraint programming allows for creating a model, that describes the problem using constraints. The possibility of using describing the expected solution, just by specifying the constraints placed on it, as well as use of the global constraints, allows for very succint and short descriptions of even very hard problems. This spirit was captured in words of E.~Freuder~\cite{Freuder1997}: 

\begin{displayquote}
``Constraint programming represents one of the closest approaches computer science has yet made to the Holy Grail of programming: the user states the problem, the computer solves it''
\end{displayquote}

\subsection*{Summary}

In case of cutting and packing problems, due to an extensive diversity of applications and requirements, the solutions presented in literature employ a variety of approaches -- from exact solutions, through approximation algorithms and heuristics, to search algorithms, as well as the approach presented also in this work: constraint programming.

Constraint programming is in principle a technique approaching the global optimum, through considering the whole space of feasible solutions, therefore it is possible, that in the case of NP-hard problems it might require exceedingly high computational resources or time, which in practice leads to imposing a time limit, after which a partial solution will be returned -- feasible, but not necessarily optimal.

\section{Solution based on Constraint Programming}
\label{ch:solution}

The main goal of the article was to describe the container packing problem using constraint programming. In this section there will be presented a description of the structure of the CP model.

\subsection{Assumptions}

Cutting and packing is a large family of problems (see \ref{sec:typology}), therefore it's imperative to specify assumptions about the exact problem for which the model will be presented.

The problem that is considered in this article is a three-dimensional, orthogonal packing, which means, that both small items and containers are rectangular cuboids (all faces are rectangles). Additionally, there is an assumption on the composition of the items: they are weakly heterogeneous -- they can be divided in a small number (relative to the overall cardinality of items) of classes of objects of the same properties. For each of the classes, there are defined axes along which the items can be rotated (at least one). The objective is to minimize the leftover payload (equivalently: maximize the packed payload).

To sum up, the problem considered in the following model falls into the \textsc{Single Large Object Placement Problem} (SLOPP) category of W\"ascher's typology, or \texttt{3/B/O/R} in Dyckhoff's classification.

\subsection{Input parameters}

The first group of input parameters consists of the container description. Container dimensions are sometimes denoted for brevity as $D_u$, where $u$ is the dimension.

\begin{center}
    \begin{tabular}{ l l }
    $W \in \NN$ & width \\
    $L \in \NN$ & length \\
    $H \in \NN$ & height \\
    \end{tabular}
\end{center}

The second group of parameters describes the payload (small items): 
\begin{center}
    \begin{tabular}{ l p{5cm} }
    $C \in \NN$ & number of classes \\
    $b_i \in \NN$ & cardinality of class $i$ \\
    $d_{i,u} \in \NN$ & size of the items of class $i$ along dimension $u$ \\
    $v_{i,u} \in \{0,1\}$ & if the object of class $i$ can be placed vertically along dimension $u$ \\
    \end{tabular}
\end{center}

\subsection{Decision variables}

Let $N$ be the number of packed items. State of item $i \in \{1, \ldots, N \}$ is described using the following decision variables:

\begin{center}
    \begin{tabular}{ l p{5cm} }
        $s_i \in \{ 1,\ldots, C \}$ & class of the item $i$, \\
        $p_{i,u} \in \NN $ & position of the item $i$ along dimension~$u$, \\
        $x_i \in \{0, 1\}$ & if the item $i$ is selected for packing, \\
        $r_i \in \mathbb{S}_3$ & item rotation. \\
    \end{tabular}
\end{center}

To simplify further definitions, a few notations were introduced: $\textsc{pos}$ describing the position of the item, taking into account rotation, and analogously $\textsc{size}$ describing the size of the item.

\[ \textsc{pos}(i, u) = p_{i, r_i(u)} \qquad \textsc{size}(i, u) = d_{s_i, r_i(u)} \]

\subsection{Objective function}

In literature, there were considered 3 types of objective functions for SLOPP:

\begin{itemize}
    \item packing density maximization (payload volume over bounding-box volume) (\cite{Loh1992}),
    \item container utilization maximization / leftover payload minimization (payload volume over container volume) (\cite{Ngoi1994, Bischoff1995, Bischoff1995a, Bortfeldt1998, Eley2002, Jin2004}),
    \item leftover items minimization / packed item maximization (ibid. and~\cite{Loh1992}).
\end{itemize}

The first of the presented objective functions can be used in cases where the volume of the payload is drastically smaller than the container, otherwise, it mostly reduces to the second function.

In the majority of the surveyed literature, the considered objective functions are the latter two -- maximization of container volume utilization and/or minimization of the number of leftover items. Under the assumption, that the value of the item is proportional to its volume, then the container utilization should be a more appropriate measure -- leaving many smaller items would leave less value behind, compared to leaving a few bigger ones.

Having considered the above points, the objective function selected for this model was the minimization of the volume of leftover payload:

\[ \sum_{i~\in \{1,\ldots,N\}} (1 - x_i) \cdot \prod_{u \in \{1, 2, 3\}} \textsc{size}(i, u) \]

\subsection{Constraints}

To simplify the notation, let $\mathfrak{S}$ be a set of classes $\{1, \ldots, C\}$, $\mathfrak{D}$ a set of dimensions, and $\mathfrak{B}$ a set of all items $\{1,\ldots,N\}$.

The first constraint class describes the constraint of number of items in each of the item classes:

\[
    |\{ i \in \mathfrak{B}:  s_i = k \}| = b_k 
        \qquad 
    \forall k \in \mathfrak{S}
\]

The above set of constraint can be rewritten using global constraint $\textsc{GlobalCardinality}(x, y, z)$, which states, that number of appearances of element $y_i$ in the array $x$ is equal to $z_i$.

\[
    \textsc{GlobalCardinality}(
        [s_i: i \in \mathfrak{B}], 
        [k: k \in \mathfrak{S}],
        [b_k: k \in \mathfrak{S}]
    ) 
\]

To reduce the number of combinations, which have isomorphic solutions (numbering of the items inside the model doesn't matter), a symmetry breaking constraint was introduced:

\[ i~< j \implies s_i \leq s_j \qquad \forall i,j \in \mathfrak{B} \]

\noindent which is realized by the global constraint \textsc{Increasing}.

The next set of constraints describes full enclosure of the items in the container, which means, that in each dimension the start of item is at least 0, and the end is at most the size of the container along that dimension:

\[
    0 \leq \textsc{pos}(i, u) 
        \land
    \textsc{pos}(i, u) + \textsc{size}(i, u) \leq D_u
        \qquad
    \forall u \in \mathfrak{D}
\]

Next set of constraints ensures, that items are properly placed, which means that the geometrical constraint (non-overlapping) is satisfied. Two rectangular cuboids don't intersect, if there exists an axis, along which the projections of the shapes don't overlap. In case of orthogonal packing, this requires only checking the primary dimensions, and the constraint for the two packed items $i,j$ can be specified as:

\begin{equation*}
\begin{aligned}
    \exists u \in \mathfrak{D}, & 
        ((\textsc{pos}(i, u) + \textsc{size}(i,u) \leq \textsc{pos}(j, u)) \lor \\
        &
        (\textsc{pos}(j, u) + \textsc{size}(j,u) \leq \textsc{pos}(i, u))) \\
\end{aligned}
\end{equation*}

Similarly as in the case of class cardinality, the constraint above can be rewritten using a global constraint, which will significantly simplify the model, and allow the solver to use specialized algorithms (cf. section \ref{sec:wiezy-globalne}). In that case, the \textsc{diffn} constraint was used\footnotemark.

\footnotetext{An implicit assumption was taken, that the items of size 0 along any dimension can be packed anywhere -- cf. Globbal Constraint Catalog, \texttt{diffn} \cite{GCC2010} In MiniZinc implementation this constraint is realized by \texttt{diffn\_nonstrict\_k}. }

\section{Tests and performance analysis}
\label{ch:tests}

\subsection{Test cases}

Performance analysis of the solution has been conducted on benchmark cases previously used in the literature. In the OR-Library\footnotemark repository, there are 7 test suites for the SLOPP problem -- \texttt{thpack1}, $\ldots$, \texttt{thpack7}, used i.a. in \cite{Bischoff1995}, and the \texttt{thpack8} suite, used i.a. in \cite{Loh1992}, \cite{Ngoi1994},\cite{Bischoff1995}, and \cite{Jin2004}.

\footnotetext{Source originally described in \cite{Beasley1990}, currently available at: \url{http://people.brunel.ac.uk/~mastjjb/jeb/info.html}}

\subsubsection{Input data format}

The first line of the input file describes the number $P$ of test cases included in the suite. Following that, provided is a description for each case $p \in {1, \ldots, P}$:
\vspace{1em}
\begin{alltt}
p [seed\textsubscript{p}]
L\textsubscript{p} W\textsubscript{p} H\textsubscript{p}
n\textsubscript{p}
1 l\textsubscript{p,1} vl\textsubscript{p,1} w\textsubscript{p,1} vw\textsubscript{p,1} h\textsubscript{p,1} vh\textsubscript{p,1} c\textsubscript{p,1}
...
n l\textsubscript{p,n} vl\textsubscript{p,n} w\textsubscript{p,n} vw\textsubscript{p,n} h\textsubscript{p,n} vh\textsubscript{p,n}
\end{alltt} 
\vspace{1em}
\noindent where \texttt{p} is an index of an instance, \texttt{seed} is an optional element describing the seed for PRNG used to generate instances (suites \texttt{thpack1-7}, \texttt{L\textsubscript{p}, W\textsubscript{p}, H\textsubscript{p}} are dimensions of the container, \texttt{n} is number of item classes, \texttt{l\textsubscript{p,i}, w\textsubscript{p,i}, h\textsubscript{p,i}} are dimensions of item from i-th class, \texttt{vl\textsubscript{p,i}, vw\textsubscript{p,i}, vh\textsubscript{p,i}} are binary variables stating if the item can be placed vertically in that dimension, and \texttt{c\textsubscript{p,i}} is the number of items in i-th class.

The test files (placed in the \texttt{data/thpack} directory) were converted to format accepted by the MiniZinc environment using the \texttt{data/util/convert.py} script. Converted instances were placed in the \texttt{data/converted} directory.

\subsubsection{Instance characterization}

Instances \texttt{thpack1-7} describe packing of a 20ft containers coming from 1C\footnotemark series of ISO 668 standard\cite{ISO668}. Instances in each suite present a different heterogeneity of the payload, e.g. \texttt{thpack1} -- 3 classes przedmiotów, \texttt{thpack4} -- 10 classes, \texttt{thpack7} -- 20 classes; number of items to pack is diversified across the instances within the suite. A summary of the characterization is presented in the table \ref{tab:char-zestawow}.

\footnotetext{Internal dimensions of $5867\text{mm} \times 2330\text{mm} \times 2197\text{mm}$, rounded to the nearest integer -- $587\text{cm} \times 233\text{cm} \times 220\text{cm}$.}

\begin{table}[!ht]
    \centering
    \rowcolors{2}{white}{black!10} 
    \begin{tabular}{ c r r }
        \toprule
        \multicolumn{1}{c}{Suite} & 
        \multicolumn{1}{c}{Classes} &
        \multicolumn{1}{c}{Items} \\ 
        \midrule
        thpack1 &  3 & 69 - 476 \\
        thpack2 &  5 & 81 - 266 \\
        thpack3 &  8 & 80 - 232 \\
        thpack4 & 10 & 75 - 233 \\
        thpack5 & 12 & 84 - 218 \\
        thpack6 & 15 & 85 - 203 \\
        thpack7 & 20 & 90 - 172 \\ 
        \bottomrule
    \end{tabular}
    \caption{Characterization of the \texttt{thpack1-7} suites.}
    \label{tab:char-zestawow}
\end{table}

The \texttt{thpack8} suite contains very diversified instances, which contrary to \texttt{thpack1-7} suites, present a payload of a volume much smaller than the volume of the container, what is caused by a different objective function in \cite{Loh1992}: maximization of packing density, i.e. volume of packed payload over volume of smallest enclosing cuboid. Results presented in \cite{Ngoi1994} employ an objective analogous to the approach used in this work, namely volume utilization (volume of the payload over volume of the container).

Instances in the \texttt{thpack8} suite describe packing of 100-250 objects, divided into 6-10 classes. In all of the test cases, each class has only one vertical dimension, which means, that it's a fragile payload (``This side up!''). Container sizes across instances are very diversified, ranging from $ 3000 \times 2000 \times  900 $ to $ 6000 \times 2800 \times 1400 $. Summary of the characterization is presented in the table \ref{tab:char-thpack8}.

\begin{table}[!ht]
    \centering
    \rowcolors{2}{white}{black!10}
    \begin{tabular} { crrl }
        \toprule
        Instance & Classes & Items & Container dimensions \\
        \midrule
         1 &  7 & 100 & $ 3000 \times 2000 \times 1000 $ \\
         2 &  8 & 200 & $ 3000 \times 2000 \times 1000 $ \\
         3 &  8 & 200 & $ 4000 \times 2400 \times 1300 $ \\
         4 &  7 & 100 & $ 3000 \times 2000 \times 1100 $ \\
         5 &  6 & 120 & $ 3000 \times 2000 \times  900 $ \\
         6 &  8 & 200 & $ 3500 \times 2400 \times 1000 $ \\
         7 &  8 & 200 & $ 3500 \times 2400 \times 1300 $ \\
         8 &  6 & 130 & $ 3200 \times 2000 \times 1200 $ \\
         9 &  8 & 200 & $ 5000 \times 2400 \times 1400 $ \\
        10 &  8 & 250 & $ 5000 \times 2400 \times 1600 $ \\
        11 &  6 & 100 & $ 3000 \times 2400 \times 1000 $ \\
        12 &  6 & 120 & $ 3200 \times 2400 \times 1000 $ \\
        13 &  7 & 120 & $ 3500 \times 2000 \times 1200 $ \\
        14 &  6 & 120 & $ 3500 \times 2200 \times 1100 $ \\
        15 & 10 & 150 & $ 6000 \times 2800 \times 1400 $ \\
        \bottomrule
    \end{tabular}
    \caption{Characterization of the \texttt{thpack8} suite.}
    \label{tab:char-thpack8}
\end{table}

\subsection{Execution environment}

All tests were conducted on the Otryt server from the Department of Fundamentals of Computer Science at the Faculty of Fundamental Problems of Technology, WUST. The server is equipped with 
Intel Xeon E7-4850 2 GHz and 256 GB RAM memory, and is working under the Debian operating system. The server allows the use of up to 80 cores, but due to a shared environment with other students and faculty researchers, only 32 cores were used. The use of the computing server allowed a high level of parallelism for the solver and therefore enabled using shorter timeouts.

The solver used for the evaluation was \texttt{fzn-or-tools} coming from Google OR-Tools\footnote{Project website: \url{https://developers.google.com/optimization}}, due to easy access, open source, and high performance, that was proven by gold medals in three categories in MiniZinc Challange 2020 \footnote{Competition website: \url{https://www.minizinc.org/challenge2020/results2020.html}}\footnote{Exceptional results were also confirmed in 2021, 2022, and 2023 editions}, a competition for constraint-programming solvers, conducted annually since 2008\cite{Minizinc2014}.

\subsection{Test results}

The testing procedure consisted of running the model in the intermediate solution mode (\texttt{-a} flag), that is, feasible, but not necessarily optimal solutions. According to the recommendation of the MiniZinc and OR-Tools docs, the \texttt{-f} flag was set, to allow the solver to ignore the hints and pick its own search strategy. To allow the analysis of the progress in time, the \texttt{-s} flag was enabled, which displays information on the running time and number of constraint propagations. The tests were run using the command:

\vspace{1em}
\begin{verbatim}
minizinc --solver ortools -a -f -s -p 32 --time-limit T model.mzn DATA
\end{verbatim}
\vspace{1em}

\noindent where \texttt{T} is a time limit in milliseconds, and \texttt{DATA} is a path to the input file of the instance.

\begin{figure}[!ht]
    \centering
    \includegraphics[width=0.95\textwidth]{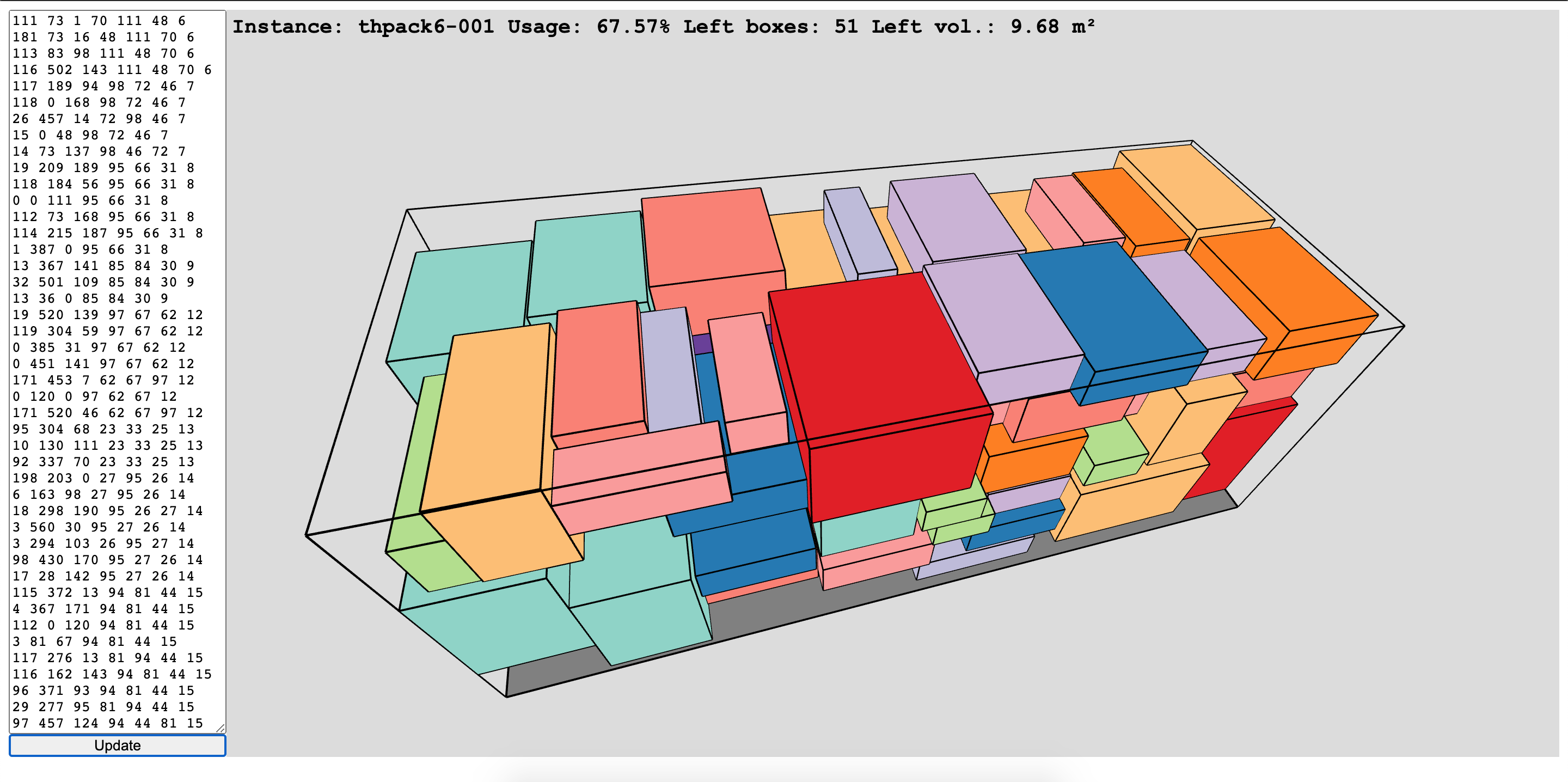}
    \caption{Example of usage of the visualization tool.}
    \label{fig:sample-result}
\end{figure}

To facilitate visual inspection of the results, a web application was created -- it's placed in the \texttt{model/visualization}. An example of usage of the app is presented in fig. \ref{fig:sample-result},

\subsubsection{Manual test}
\label{sub:manual-test}

\begin{figure}[!ht]
    \centering
    \includegraphics[width=0.95\textwidth]{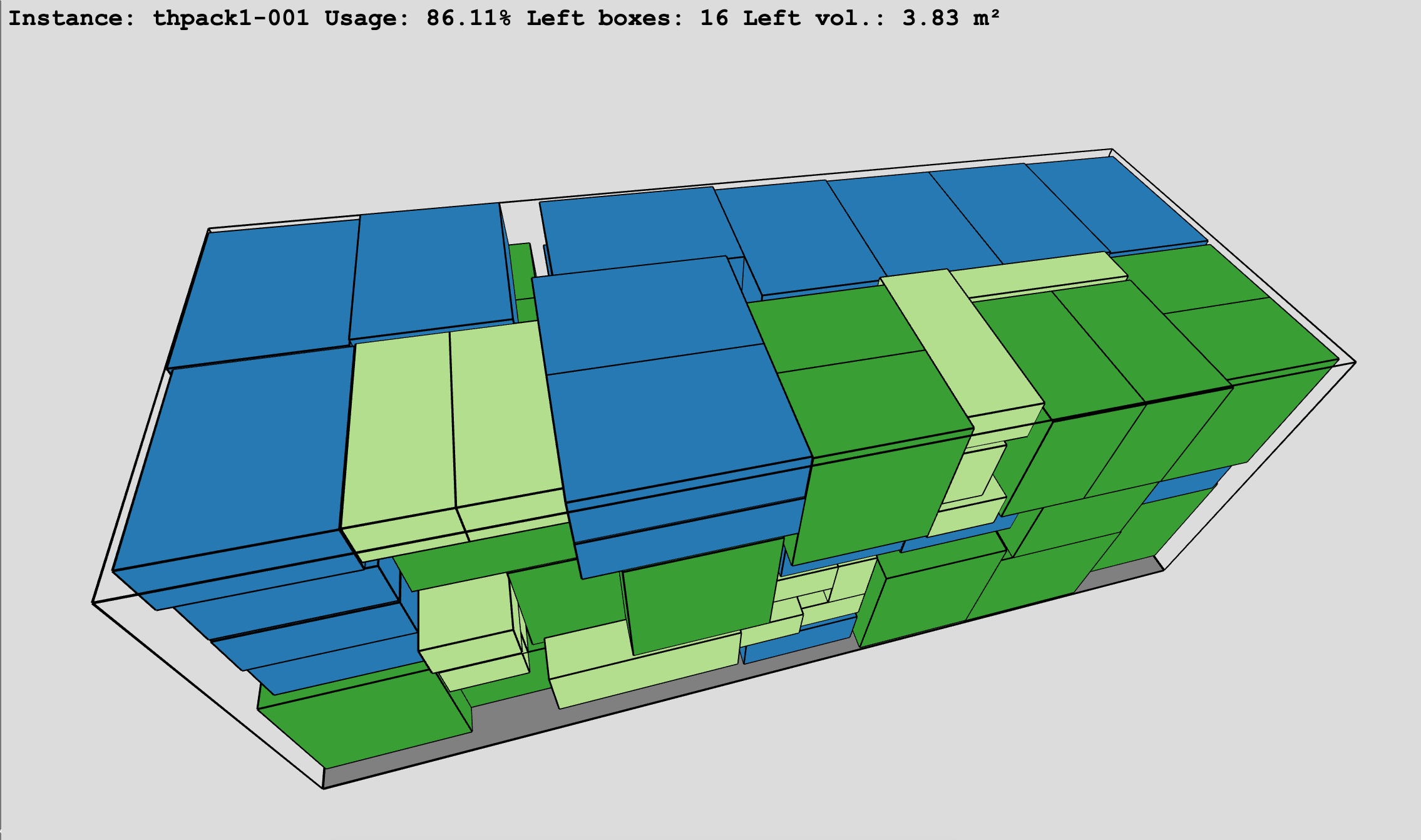}
    \caption{Manual test: visualisation of the best solution}
    \label{fig:manual-vis}
\end{figure}

To check the model performance with a higher time limit, a manual test with instance \texttt{thpack1\_001} was performed. The planned timeout was 4 hours with 32 cores used. A visualisation of the best achieved solution is presented in Fig. \ref{fig:manual-vis}

\begin{figure}[!ht]
    \centering
    \includegraphics[width=0.7\textwidth]{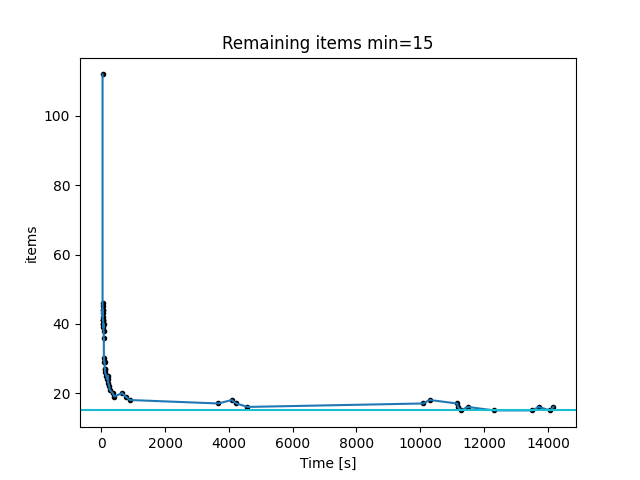}
    \caption{Manual test: remaining items}
    \label{fig:manual-boxes}
\end{figure}

\begin{figure}[!ht]
    \centering
    \includegraphics[width=0.7\textwidth]{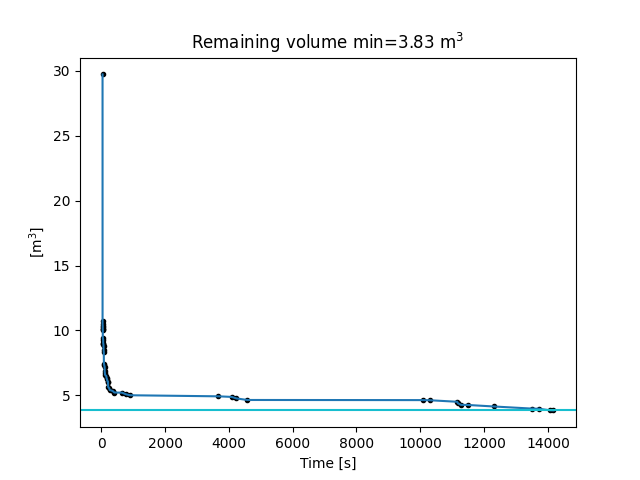}
    \caption{Manual test: remaining volume}
    \label{fig:manual-vol}
\end{figure}

\begin{figure}[!ht]
    \centering
    \includegraphics[width=0.7\textwidth]{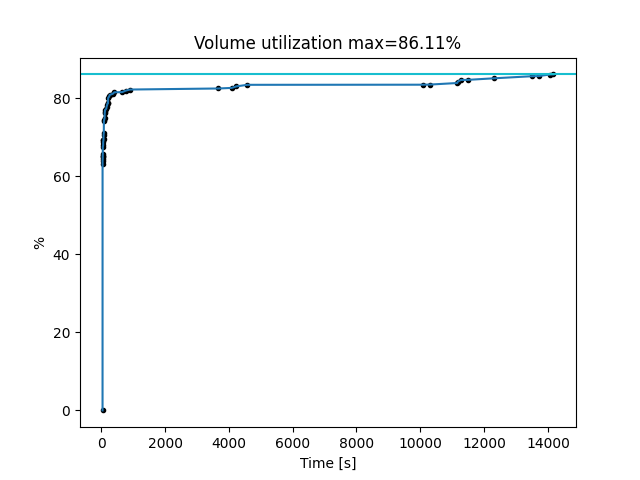}
    \caption{Manual test: volume utilization}
    \label{fig:manual-usage}
\end{figure}

Based on the data collected in the experiment, the figures \ref{fig:manual-boxes}, \ref{fig:manual-vol}, \ref{fig:manual-usage} present the values of metrics, respectively: number of remaining items, the volume of remaining payload, and volume utilization. Raw output from the data is placed in the \texttt{results/rawdata/manual} directory. The figures were plotted using the \texttt{results/script/manual\_results.py} script. 

The data gathered in the experiment show, that progress after the first 30 minutes ($82.23\%$ after 899s) is very similar to the progress achieved in 90 minutes ($83.44\%$ after 4582s), but the final results after 4 hours presented a significant improvement ($86.11\%$ after 14174s).

\subsubsection{OR-Library thpack1-7 instances}
\label{sub:test-thpack1-7}

Based on the manual testing (cf. section \ref{sub:manual-test}) and due to the limited time of the project, the time limit for the instances from \texttt{thpack1,...,thpack7} suites was set to 30 minutes, with 32 cores, and first 10 cases from each suite were selected.

The collected results were placed in the \texttt{results/rawdata/thpack1} - \texttt{thpack7} directories. In the plots presented are the gathered metrics of intermediate solutions: volume utilization, remaining items, and remaining volume. Example plot for \texttt{thpack1} suite was presented in fig. \ref{fig:thpack1_inline}. The figures for all of the suites are presented in the appendix \ref{wykresy}. In the plots, the intermediate results for each instance are presented in the differently-coloured series. A summary of the results, comparing them to prior literature, is presented in the table \ref{tab:podsumowanie}.

\begin{figure}[!ht]
\subfloat[Volume utilization]{
    \label{fig:thpack1_usage_inline} 
    \includegraphics[width=0.32\textwidth]{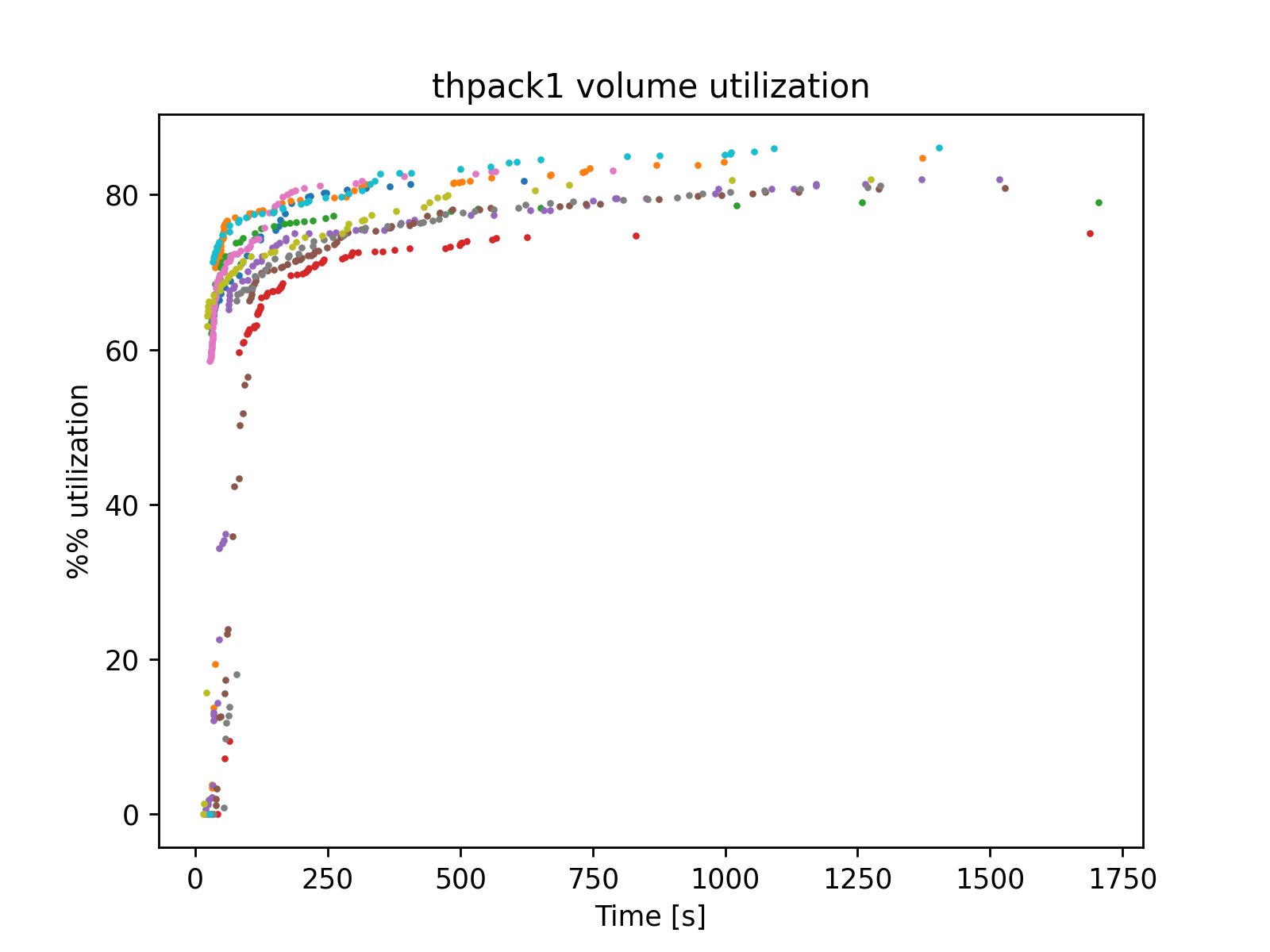}
}
\subfloat[Remaining items]{
    \label{fig:thpack1_boxes_inline}    
    \includegraphics[width=0.32\textwidth]{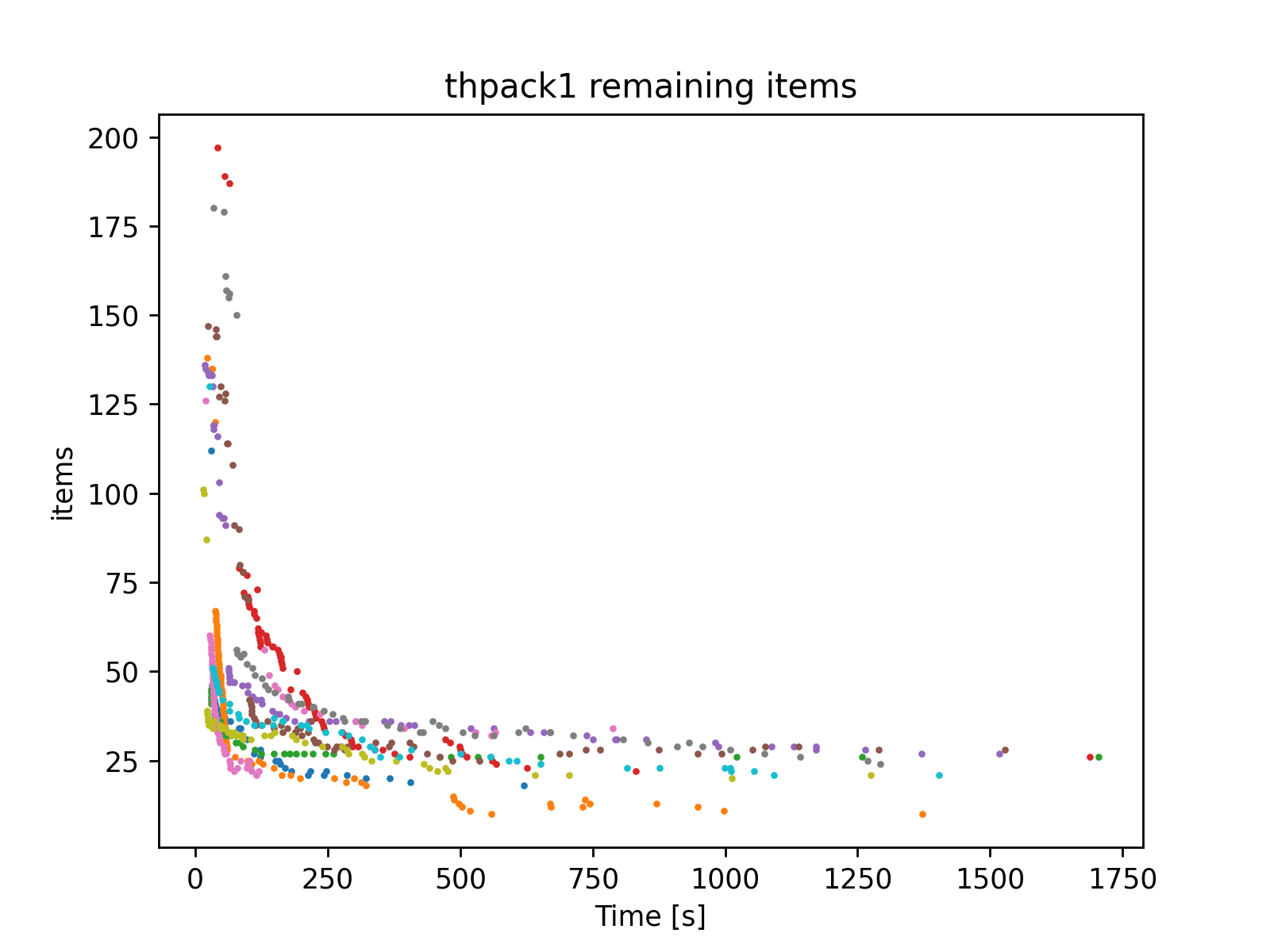}
}
\subfloat[Remaining payload]{
    \label{fig:thpack1_vol_inline}    
    \includegraphics[width=0.32\textwidth]{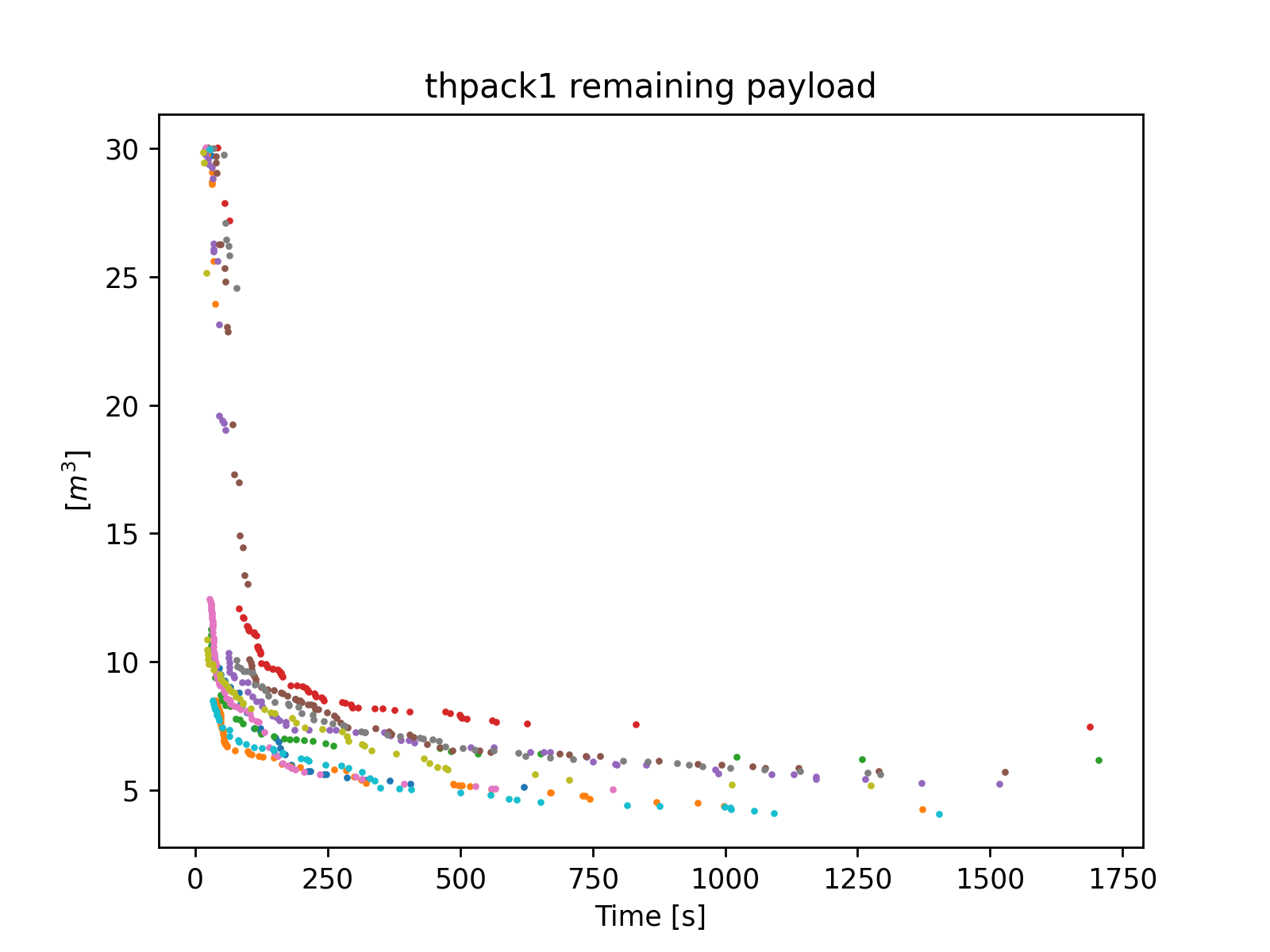}
}
\caption{Plots of results for \texttt{thpack1} test cases}
\label{fig:thpack1_inline}
\end{figure}

\begin{table}[!ht]
    \centering
    \rowcolors{3}{white}{black!10}
    \begin{tabular}{ c ccc p{10pt} ccc p{10pt} ccc }
        \toprule
        \multirow{2}{*}{Case} & 
            \multicolumn{3}{ c }{Volume utilization $[\%]$} &  ~&
            \multicolumn{3}{ c }{Remaining items} & ~&
            \multicolumn{3}{ c }{Remaining payload $[m^3]$} \\
        \cmidrule{2-4}
        \cmidrule{6-8}
        \cmidrule{10-12}
        & min & avg & max &~
        & min & avg & max &~
        & min & avg & max \\ 
        \midrule
        thpack1 &
            74.98 & 81.58 & 86.09 & ~&
            10 & 23.50 & 34 & ~&
            4.06 & 5.38 & 7.48 \\
        thpack2 &
            77.74 & 82.23 & 85.22 & ~&
            6 & 27.10 & 48 & ~&
            4.29 & 5.16 & 6.39 \\
        thpack3 &
            76.72 & 81.31 & 84.73 & ~&
            12 & 23.70 & 45 & ~&
            4.49 & 5.49 & 6.69 \\
        thpack4 &
            76.34 & 81.07 & 83.94 & ~&
            13 & 22.40 & 30 & ~&
            4.11 & 5.47 & 7.05 \\
        thpack5 &
            78.28 & 81.55 & 83.87 & ~&
            13 & 21.50 & 28 & ~&
            4.46 & 5.37 & 6.37 \\
        thpack6 &
            77.10 & 80.95 & 83.87 & ~&
            16 & 21.80 & 31 & ~&
            4.22 & 5.50 & 6.34 \\
        thpack7 &
            80.06 & 82.94 & 85.55 & ~&
            16 & 20.20 & 28 & ~&
            4.26 & 4.96 & 6.00 \\
        \bottomrule

    \end{tabular}
    \caption{Summary of results from \texttt{thpack1-7} suite}
    \label{tab:podsumowanie}
\end{table}
    
\subsubsection{OR-Library thpack8 suite}
\label{sub:res-thpack8}

Analogically to \texttt{thpack1-7} suites, tests were conducted with a time limit of 30 minutes, on 32 cores. Due to the diverse character of instances and the goal of comparison with the previous literature results, a summary of the results is presented in the table \ref{tab:results-thpack8}.

\begin{sidewaystable}[!ht]
    \rowcolors{3}{white}{black!10}
    \begin{tabular} { |r| *{7}{
            *{2}{ | >{\raggedleft\arraybackslash}p{0.7cm} }
        | } |
    }
        \hline
            Instance & 
            \multicolumn{2}{ m{1.4cm}|| }{\centering 
                Wróbel (2021)
            } &
            \multicolumn{2}{ m{1.4cm}|| }{\centering
                Jin et~al. \cite{Jin2004} (2004)
            }  &
            \multicolumn{2}{ m{1.4cm}|| }{\centering
                Eley \cite{Eley2002} (2002)
            }  &
            \multicolumn{2}{ m{1.4cm}|| }{\centering
                Bortfeldt, Gehring \cite{Bortfeldt1998} (1998)
            }  &
            \multicolumn{2}{ m{1.4cm}|| }{\centering
                Bischoff, Ratcliff \cite{Bischoff1995} (1995)
            }  &
            \multicolumn{2}{ m{1.4cm}|| }{\centering
                Bischoff et~al. \cite{Bischoff1995a} (1995)
            }  &
            \multicolumn{2}{ m{1.4cm}|| }{\centering
                Ngoi et~al. \cite{Ngoi1994} (1994)
            } \\
        \hline
            Measure 
            & \multicolumn{1}{c|}{VU} & \multicolumn{1}{c||}{LB}
            & \multicolumn{1}{c|}{VU} & \multicolumn{1}{c||}{LB}
            & \multicolumn{1}{c|}{VU} & \multicolumn{1}{c||}{LB}
            & \multicolumn{1}{c|}{VU} & \multicolumn{1}{c||}{LB}
            & \multicolumn{1}{c|}{VU} & \multicolumn{1}{c||}{LB}
            & \multicolumn{1}{c|}{VU} & \multicolumn{1}{c||}{LB}
            & \multicolumn{1}{c|}{VU} & \multicolumn{1}{c||}{LB}
            \\
        \hline
         1 & 
         62.50 &  0 &
         62.5  &  0 &
         62.5  &  0 &
         62.5  &  0 & 
         62.5  &  0 &
         62.5  &  0 &
         62.5  &  0 \\
        
         2 & 
         79.09 & 42 &
         87.7  & 28 &
         90.8  & 53 &
         96.7  & 28 & 
         90.0  & 35 &
         89.7  & 23 &
         80.73 & 54 \\
       
         3 & 
         53.43 &  0 &
         53.4  &  0 &
         53.4  &  0 &
         53.4  &  0 & 
         53.4  &  0 &
         53.4  &  0 &
         53.43 &  0 \\
        
         4 & 
         54.96 &  0 &
         55.0  &  0 &
         55.0  &  0 &
         55.0  &  0 & 
         55.0  &  0 &
         55.0  &  0 &
         54.96 &  0 \\
        
         5 & 
         77.19 &  0 &
         77.2  &  0 &
         77.2  &  0 &
         77.2  &  0 & 
         77.2  &  0 &
         77.2  &  0 &
         77.19 &  0 \\
        
         6 & 
         76.22 & 51 &
         87.4  & 28 &
         87.9  & 44 &
         96.2  & 32 & 
         83.1  & 77 &
         89.5  & 24 &
         88.72 & 48 \\
        
         7 & 
         70.12 & 29 &
         83.3  &  3 &
         84.7  &  0 &
         84.7  &  0 & 
         78.7  & 18 &
         83.9  &  1 &
         81.81 & 10 \\
        
         8 & 
         59.42 &  0 &
         59.4  &  0 &
         59.4  &  0 &
         59.4  &  0 & 
         59.4  &  0 &
         59.4  &  0 &
         59.42 &  0 \\
        
         9 & 
         61.89 &  0 &
         61.9  &  0 &
         61.9  &  0 &
         61.9  &  0 & 
         61.9  &  0 &
         61.9  &  0 &
         61.89 &  0 \\
       
        10 & 
         63.59 & 12 &
         67.3  &  0 &
         67.3  &  0 &
         67.3  &  0 & 
         67.3  &  0 &
         67.3  &  0 &
         67.29 &  0 \\
        
        11 & 
         62.16 &  0 &
         62.2  &  0 &
         62.2  &  0 &
         62.2  &  0 & 
         62.2  &  0 &
         62.2  &  0 &
         62.16 &  0 \\
        
        12 & 
         78.02 &  1 &
         78.5  &  0 &
         78.5  &  0 &
         78.5  &  0 & 
         78.5  &  0 &
         76.5  &  3 &
         78.51 &  0 \\
        
        13 & 
         79.37 &  6 &
         85.6  &  0 &
         85.6  &  0 &
         85.6  &  0 & 
         78.1  & 20 &
         82.3  &  5 &
         84.14 &  2 \\
        
        14 & 
         62.81 &  0 &
         62.8  &  0 &
         62.8  &  0 &
         62.8  &  0 &
         62.8  &  0 &
         62.8  &  0 &
         62.81 &  0 \\
        
        15 & 
         58.89 &  3 &
         59.5  &  0 &
         59.5  &  0 &
         59.5  &  0 & 
         59.5  &  0 &
         59.5  &  0 &
         59.46 &  0 \\
        \hline
        \hiderowcolors
        \textbf{TOTAL} & 
        66.64 & 144 &
        69.6  & 59 &
        69.9  & 97 &
        70.9  & 60 &
        68.6  & 150 &
        69.5  & 56 &
        69.0  & 114 \\
        \hline
        
        \multicolumn{15}{|l|}{} \\
        \multicolumn{15}{|l|}{\textbf{Notations:}} \\
        \multicolumn{15}{|l|}{
            \quad VU -- Volume Utilization,
        } \\
        \multicolumn{15}{|l|}{
            \quad LB -- Left Boxes,
        } \\
        \multicolumn{15}{|l|}{
            \quad TOTAL -- summary of the results for all instances: mean VU, sum of LB
        } \\
        \multicolumn{15}{|l|}{} \\
        \hline
    \end{tabular}
    \caption{Comparison of the results for \texttt{thpack8} suite.}
    \label{tab:results-thpack8}
\end{sidewaystable}

\subsection{Discussion}

The results of experiments in the section \ref{sub:res-thpack8} show, that the performance of the model presented in this work, under the specified resource constraints (30 minute time limit, 32 cores) has not achieved results significantly better than previous literature. 

In the case of previous works: for instances 1,3,4,5,8,9,10,11,14,15 all approaches reached optimal solution, for instance 12 all but \cite{Bischoff1995a} reached optimum, for instance 13 optimal solution was found by approaches of \cite{Jin2004}, \cite{Eley2002}, and \cite{Bortfeldt1998}, and for instance 7 only \cite{Eley2002} i~\cite{Bortfeldt1998} reached optimum. In case of instances 2 and 6 each approach has given different results.

Comparing the above results with the experiments conducted using the model proposed in this work, it's apparent, that in case of instances 1, 3, 4, 5, 8, 9, 11, 14, similarly as in the above sources the optimum was achieved. In cases 10 and 15, contrary to other solutions an optimal value was not achieved. In cases 12 and 13 optimal value was not achieved, but according to both objective functions, final solution was better than solutions from \cite{Bischoff1995a} and \cite{Bischoff1995}. 

In instance 2, proposed model has performed worse than others, with regard to volume utilization, but has left less items than solutions  \cite{Eley2002} and \cite{Ngoi1994}. In case of instance 6 the model has given significantly worse solution (almost 7 percentage points to the next solution), but has left fewer items behind than \cite{Bischoff1995}. In case 7, the achieved solution was worse than prior literature under both metrics.

Having considered the whole of the test suite, the results of the proposed model are worse than results achieved previously in literature, but in some cases it has left fewer items behind. Taking into the account results of manual test (cf. section  \ref{sub:manual-test}) increasing the time limit or number of used cores might lead to improvement of the achieved results.

\section{Summary}
\label{ch:summary}

Cutting and packing problems are ubiquitous in many areas of life and industry, therefore they appear in literature under different names. The three main groups of problems, mentioned in this work, were 3D Bin Packing (3D-BPP), pallet packing (3D Open Dimension Problem, 3D Strip Packing, 3D-ODP), and 3D Single Container Packing (3D-SLOPP)

The main goals of this work, that is: to specify a Constraint-Programming-based model for the 3D Single Container Packing problem and compare its performance with previous literature results. The initial hypothesis, that the Constraint Programming model's performance is worse than heuristic heuristic and metaheuristic algorithms, was confirmed by the results of experiments on widely used datasets. A possible explanation of such an outcome is that CP-based algorithms are searching for globally optimal solutions, which requires browsing a much bigger search space compared to local search algorithms. The results presented in this work were gathered in a shared environment, with limited computational resources, and it was necessary to impose a time limit for each of the test cases. Additionally, as was shown in the example of the single case with the extended time limit, increasing the time limit resulted in a significant improvement in the solution. 

In the future, the author is planning to continue the research on possibilities of using constraint programming in other operational research problems, as well as using hybrid search algorithms (see \cite{Wallace2007}), which may allow a combination of the global optimization of constraint programming with the performance of local search algorithms. Another avenue for research is robust optimization, where the input values or constraints might be uncertain (e.g. measurement errors, unknown or hard-to-estimate processes), applied both to constraint programming as well as cutting and packing problems.

One of the topics that wasn't thoroughly researched in this work is the topic of problem variants arising from practical contexts. In the model presented in this work, only the verticality constraint was introduced (``This side up!''), but in the literature many others were discussed: load bearing, stability criteria, and multiple drop points (Multi-drop Constraint). A rarely mentioned in the literature, but having significant industrial applications, a variant of the packing problem is the online formulation, where the assortment is not fully known at the start, and dynamic formulation, where objects are dynamically arriving and departing.

Another interesting topic, that was only mentioned in this work, is the usage of redundant constraints and symmetry-breaking constraints, which allow for a better exploitation of the known structure of the problem (e.g. weak heterogeneity of the payload), which may lead to reducing the search space by a trade-off of an increasing number of constraints and complexity of the model. An interesting question, that has arisen during the initial work on the model is how the symmetry-breaking constraints and redundant constraints impact the performance of the model and the quality of the results.

To sum up, the primary goals of the work, that is presenting and implementing the constraint-programming-based model for 3D packing, have been achieved. Despite the limited timespan of the project, the model was tested with a self-imposed time limit, and the gathered results are comparable to prior literature results. The proposed further work areas may lead to performance improvement, and allow the implementation in the practical business applications.

\backmatter

\section*{Disclaimer}

The research presented in this work was originally conducted as a part of the author's Master's Thesis under the guidance of prof. dr hab. Paweł Zieliński, in the Faculty of Fundamental Problems of Science at Wrocław University of Science and Technology.

\section*{Supplementary material}

The code and results mentioned in the work are available online in the GitHub repository: \url{https://github.com/KatJon/3d-packing}.

\FloatBarrier

\begin{appendices}
\section{Results and plots}
\label{wykresy}

Plots presented in sections \ref{sec:appendix-usage}, \ref{sec:appendix-boxes}, and \ref{sec:appendix-vol} contain results of experiments described in section \ref{sub:test-thpack1-7} for all test suites (\texttt{thpack1} - \texttt{thpack7}).

\subsection{Volume usage}
\label{sec:appendix-usage}

\begin{figure}[!ht]
    \centering
    \includegraphics[width=0.75\linewidth]{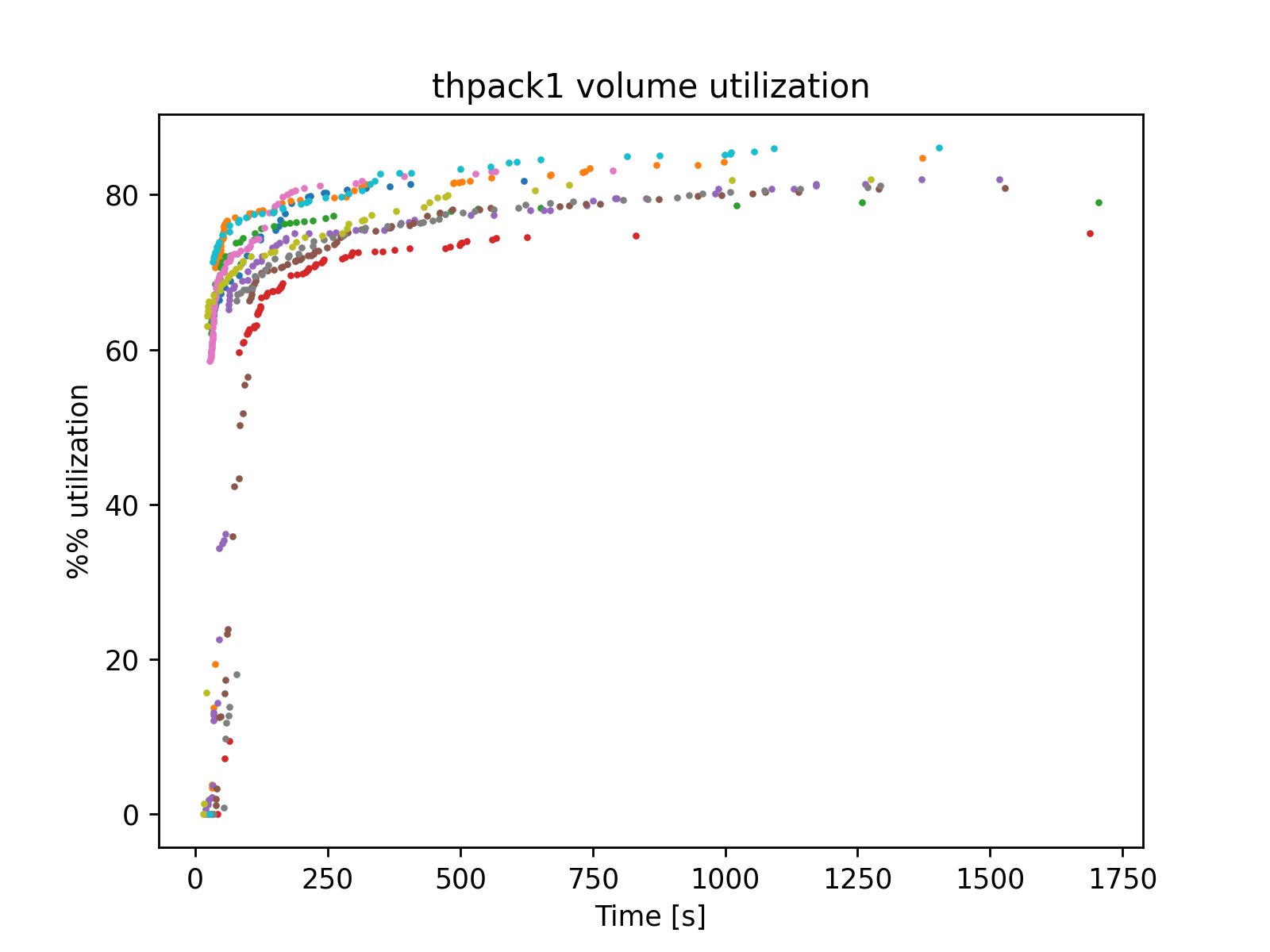}
    \caption{Volume utilization \texttt{thpack1}}
    \label{fig:thpack1_usage}    
\end{figure}

\begin{figure}[!ht]
    \centering
    \includegraphics[width=0.75\linewidth]{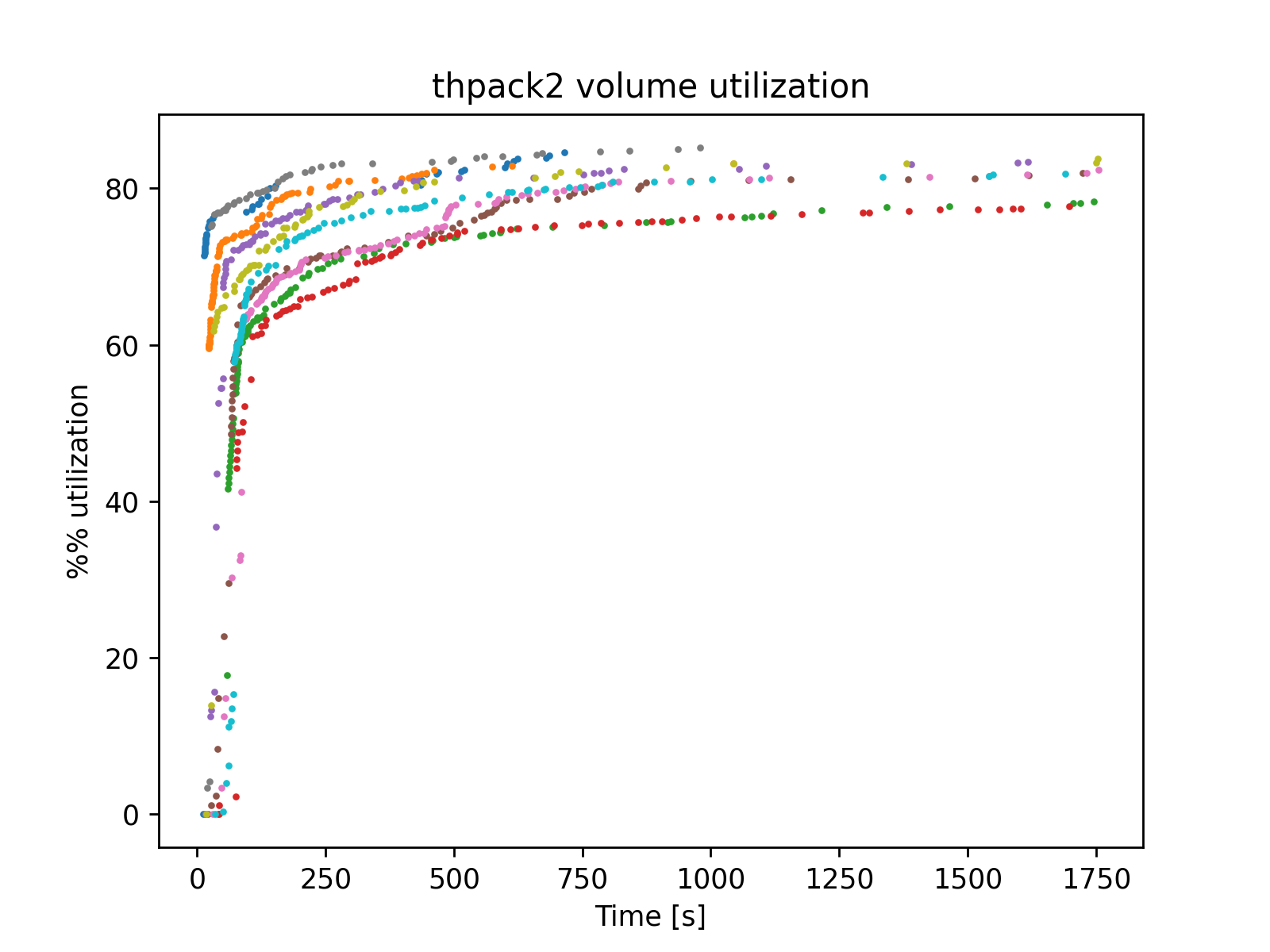}
    \caption{Volume utilization \texttt{thpack2}}
    \label{fig:thpack2_usage}    
\end{figure}

\begin{figure}[!ht]
    \centering
    \includegraphics[width=0.75\linewidth]{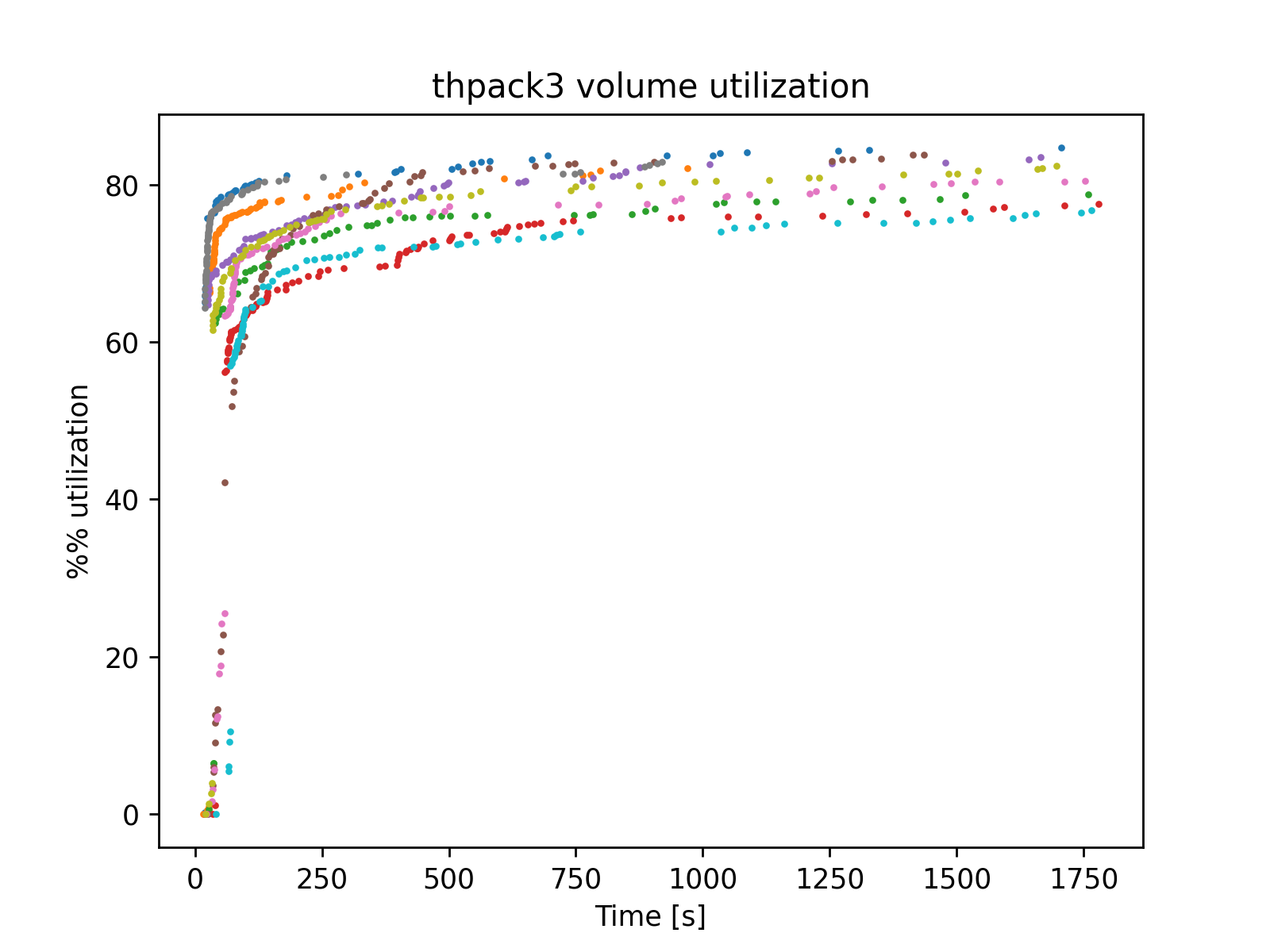}
    \caption{Volume utilization \texttt{thpack3}}
    \label{fig:thpack3_usage}    
\end{figure}

\begin{figure}[!ht]
    \centering
    \includegraphics[width=0.75\linewidth]{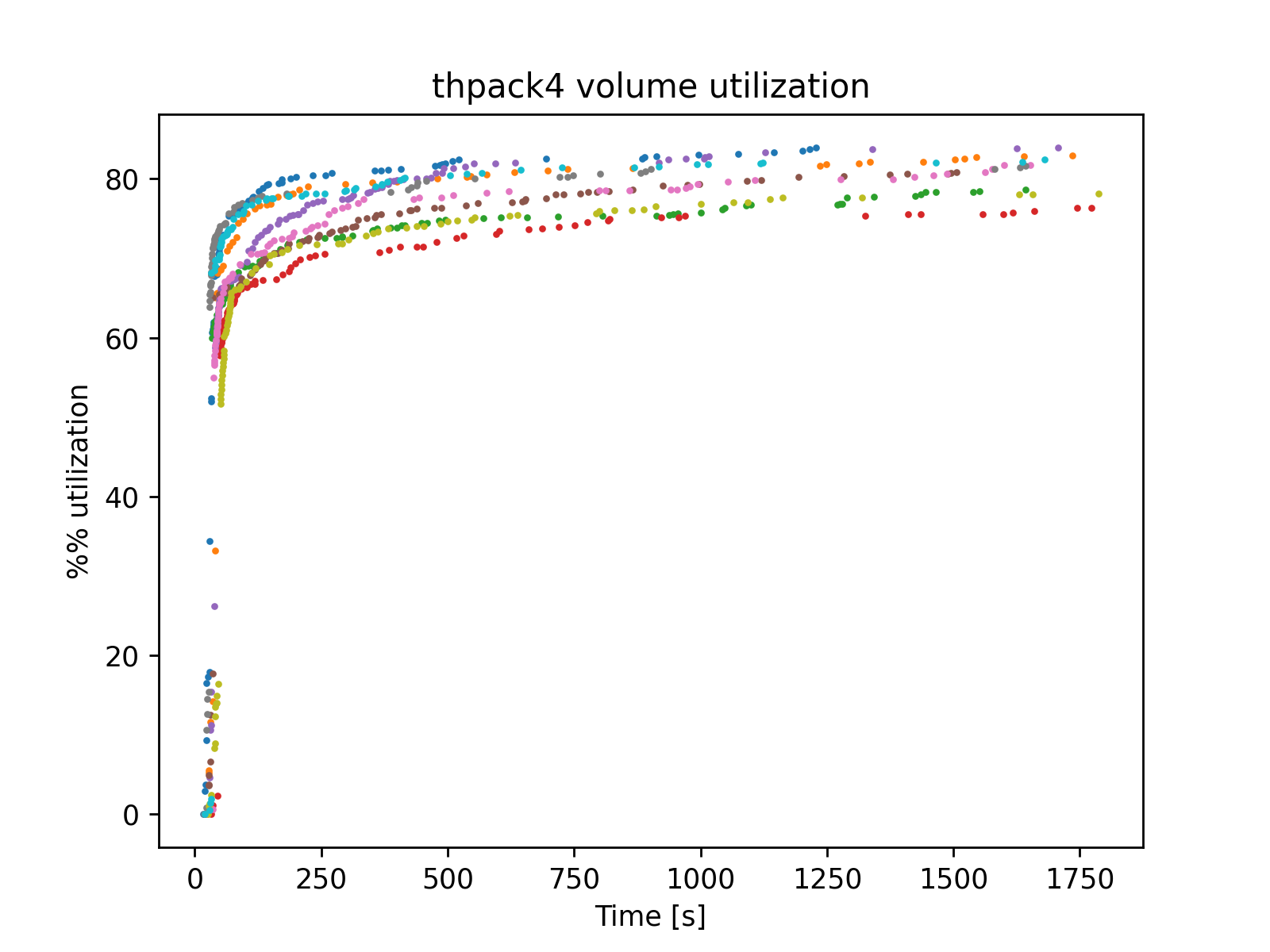}
    \caption{Volume utilization \texttt{thpack4}}
    \label{fig:thpack4_usage}    
\end{figure}

\begin{figure}[!ht]
    \centering
    \includegraphics[width=0.75\linewidth]{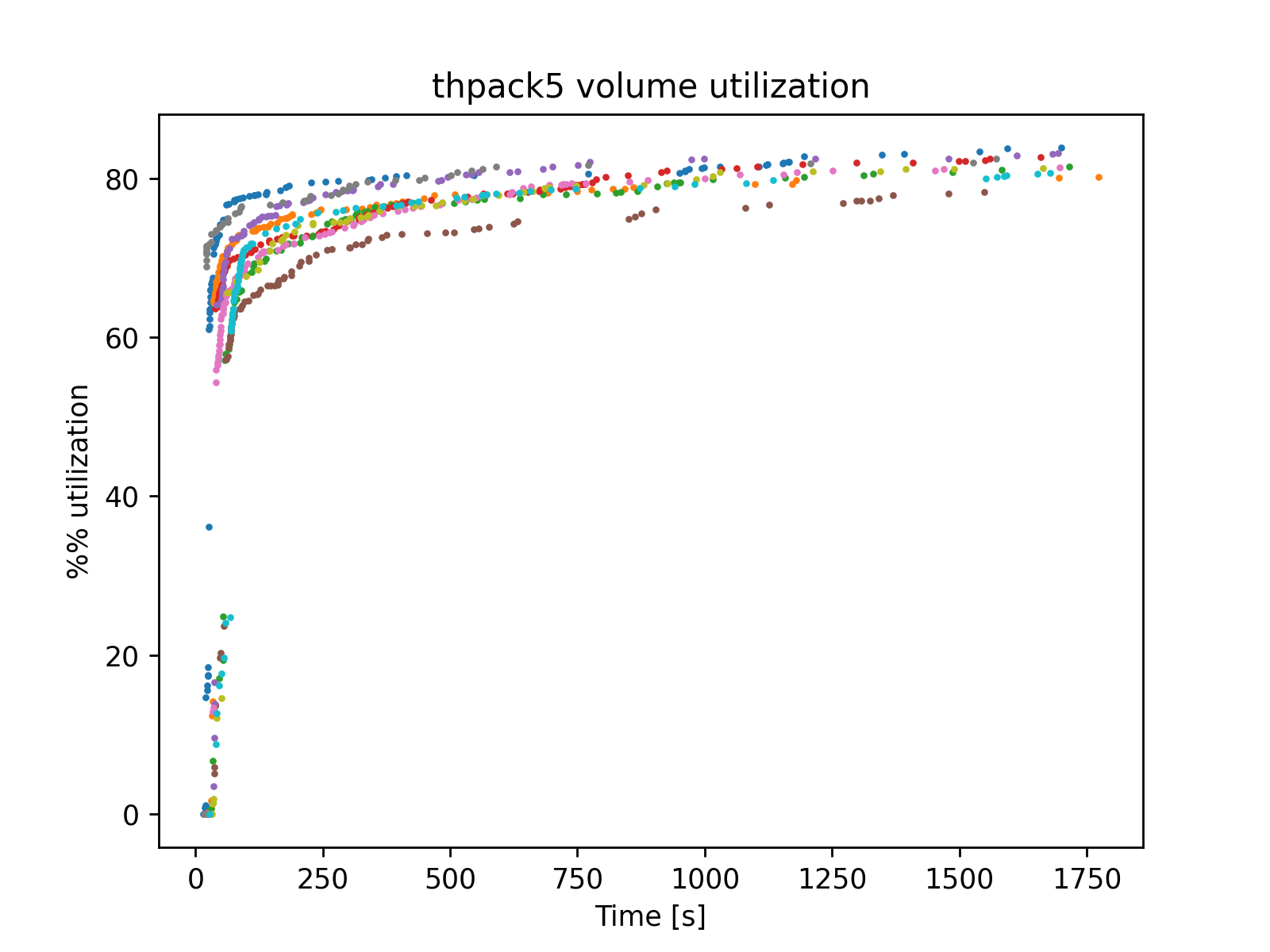}
    \caption{Volume utilization \texttt{thpack5}}
    \label{fig:thpack5_usage}    
\end{figure}

\begin{figure}[!ht]
    \centering
    \includegraphics[width=0.75\linewidth]{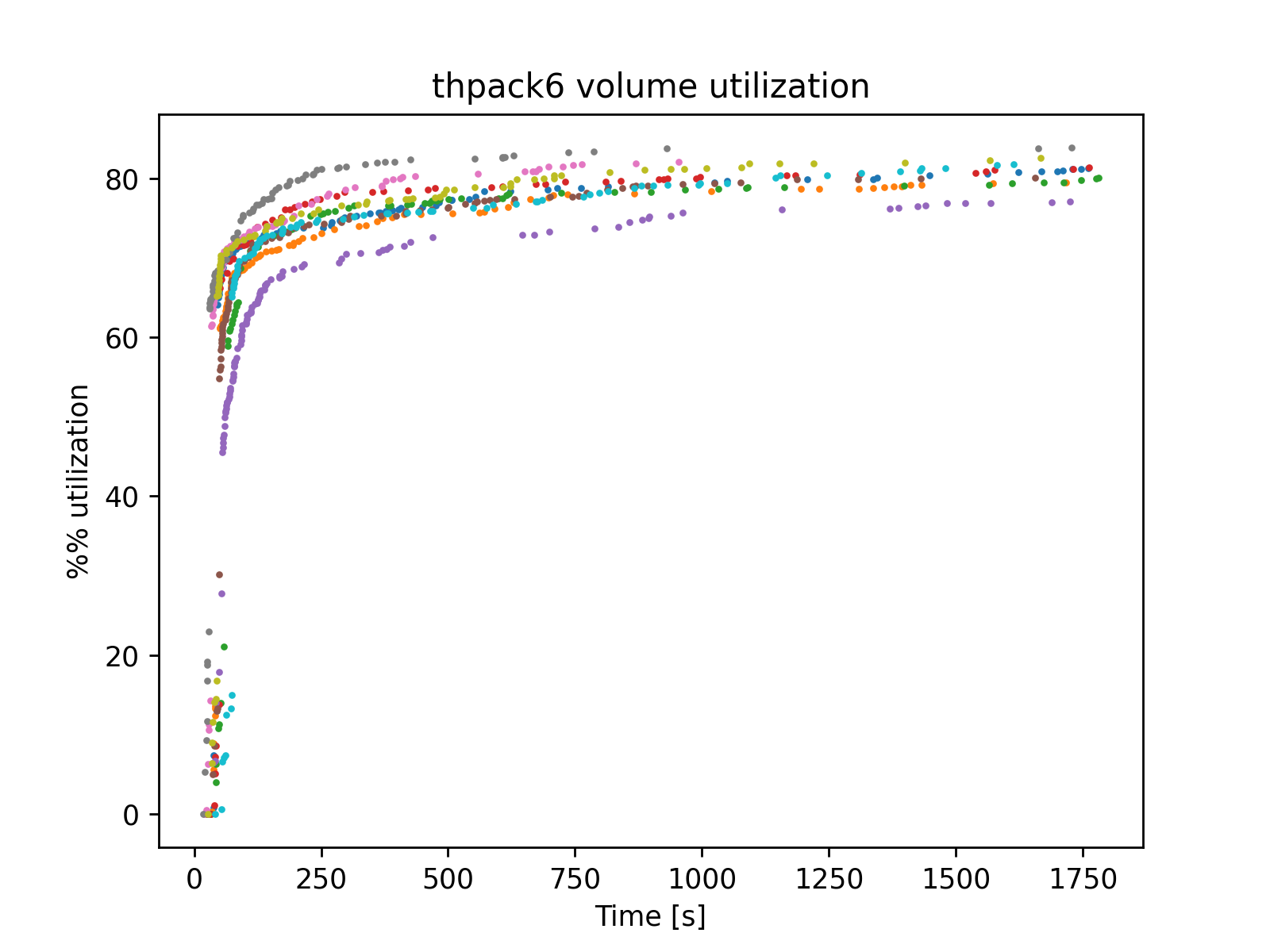}
    \caption{Volume utilization \texttt{thpack6}}
    \label{fig:thpack6_usage}    
\end{figure}

\begin{figure}[!ht]
    \centering
    \includegraphics[width=0.75\linewidth]{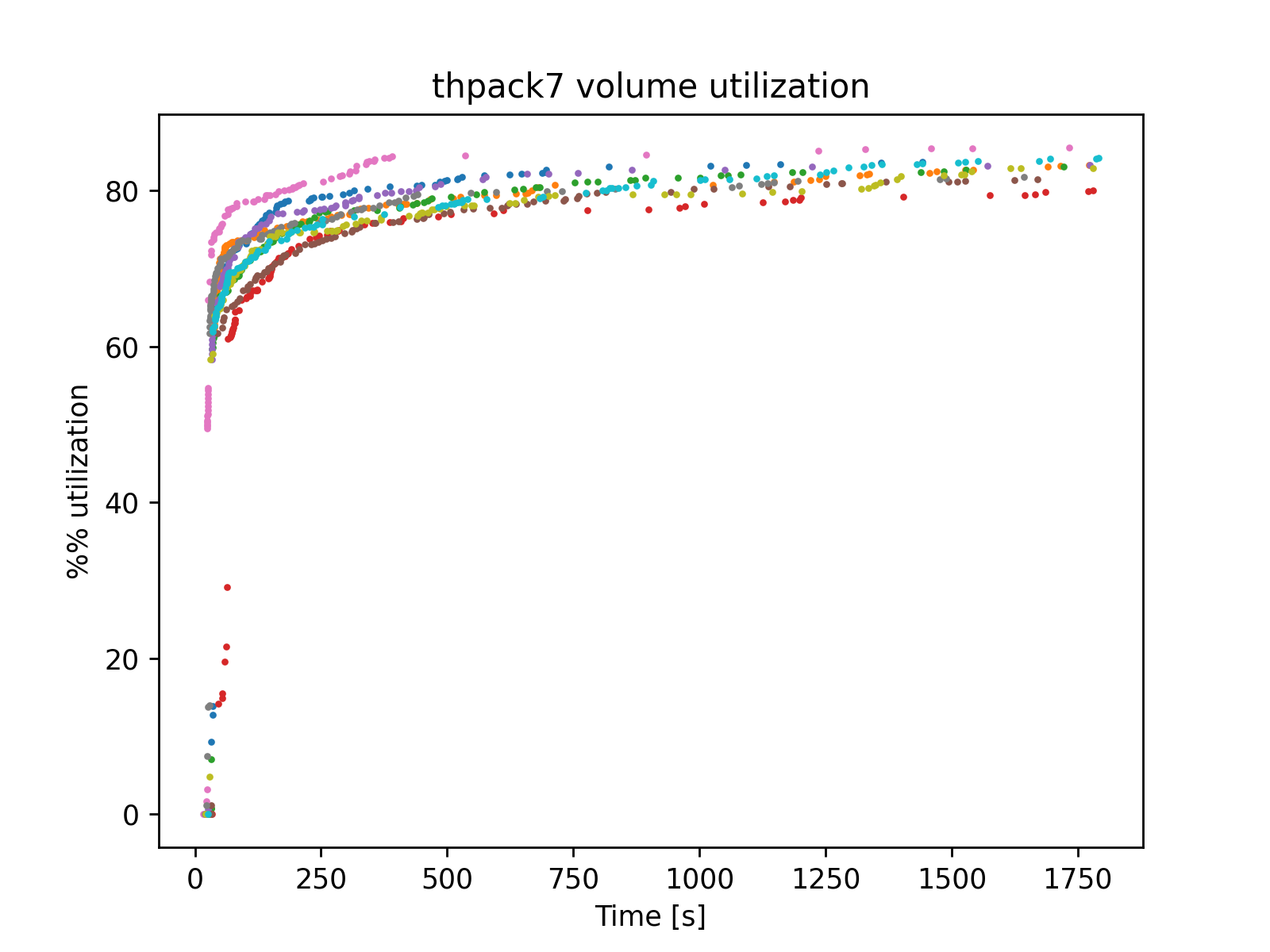}
    \caption{Volume utilization \texttt{thpack7}}
    \label{fig:thpack7_usage}    
\end{figure}

In plots \ref{fig:thpack1_usage}-\ref{fig:thpack7_usage} presented is volume utilization of partial solutions gathered from experiments in section \ref{sub:test-thpack1-7}.

\FloatBarrier

\subsection{Remaining items}
\label{sec:appendix-boxes}

\begin{figure}[!ht]
    \centering
    \includegraphics[width=0.75\linewidth]{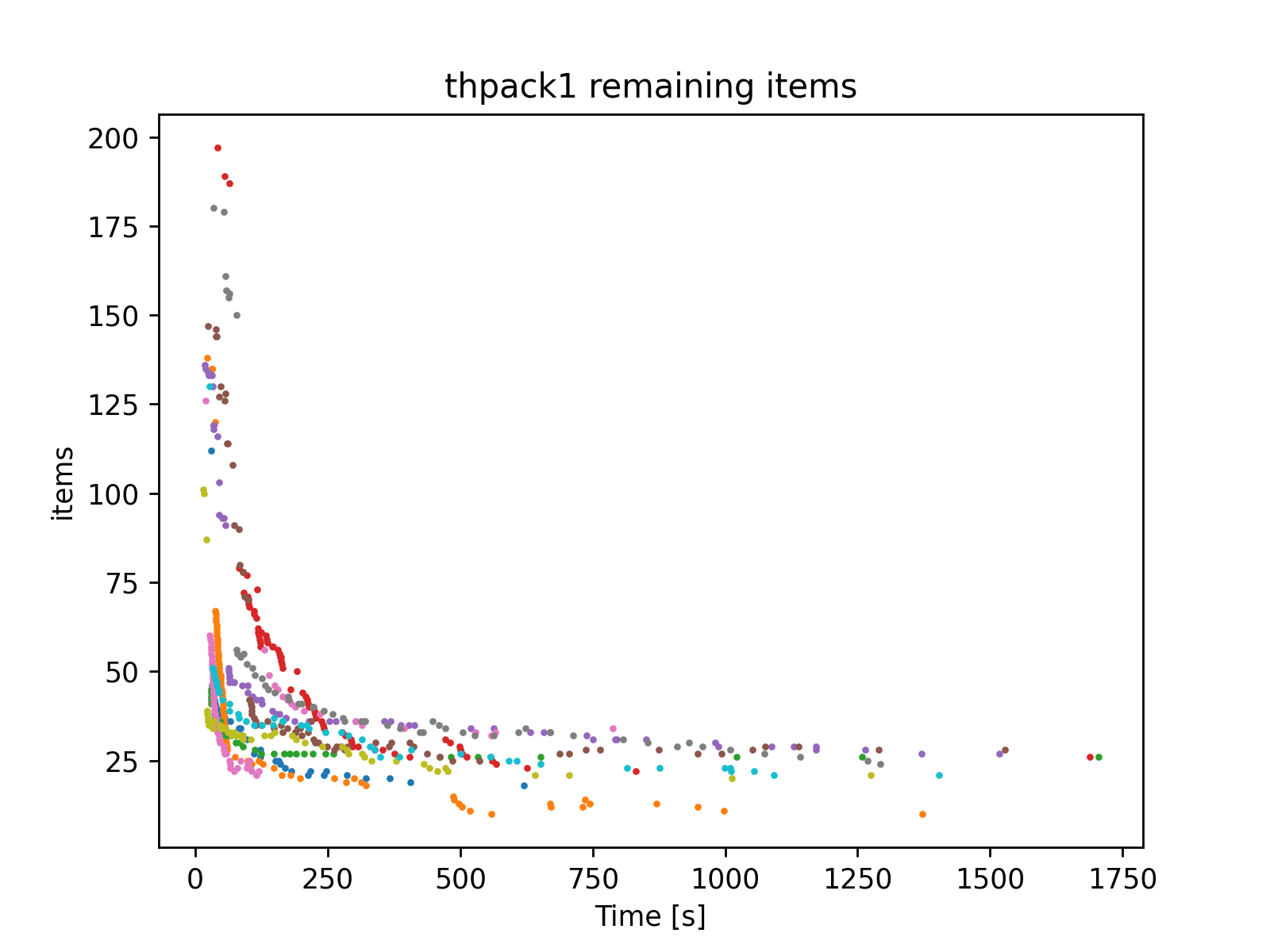}
    \caption{Remaining items \texttt{thpack1}}
    \label{fig:thpack1_boxes}    
\end{figure}

\begin{figure}
    \centering
    \includegraphics[width=0.75\linewidth]{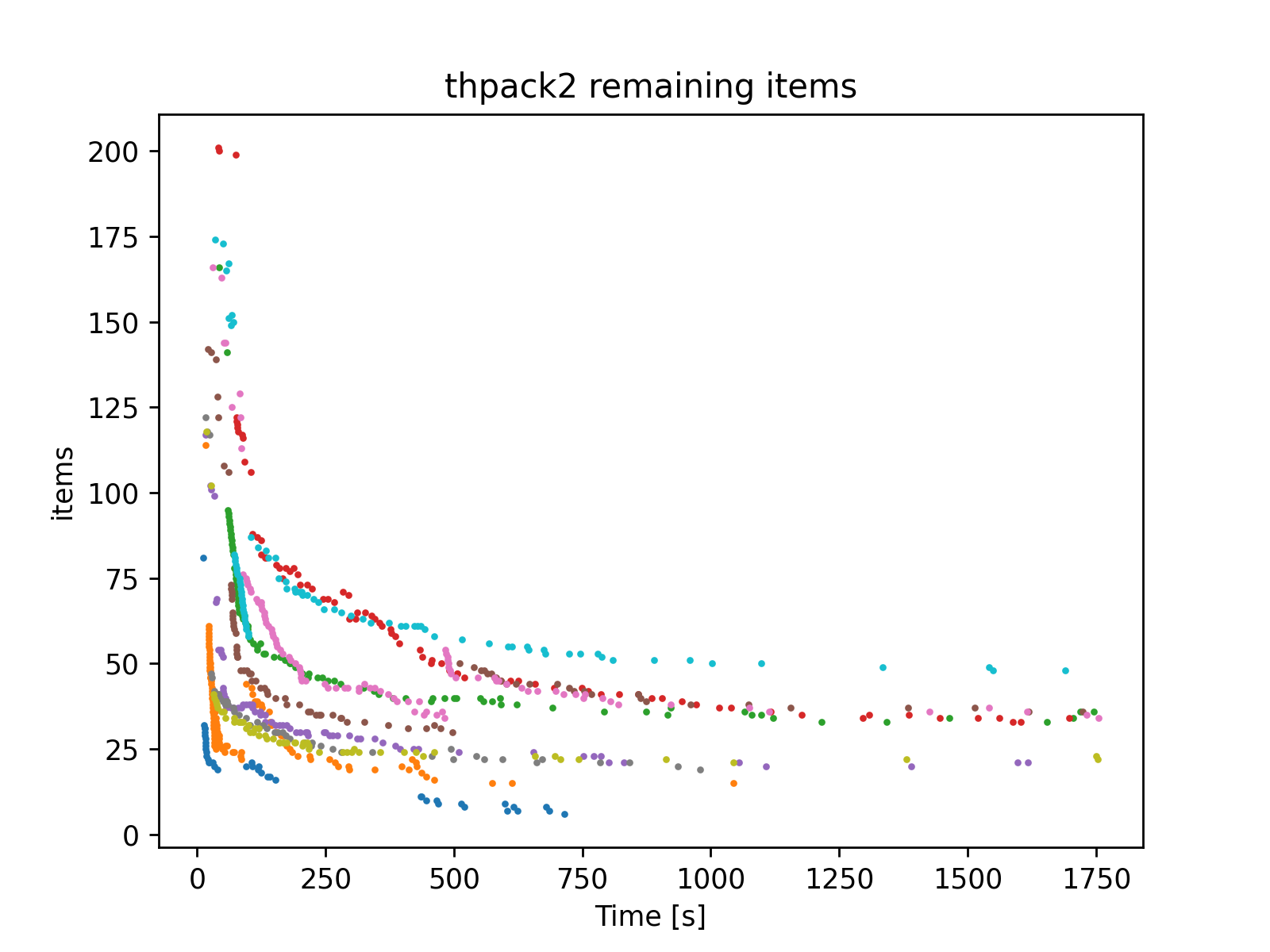}
    \caption{Remaining items \texttt{thpack2}}
    \label{fig:thpack2_boxes}    
\end{figure}

\begin{figure}[!ht]
    \centering
    \includegraphics[width=0.75\linewidth]{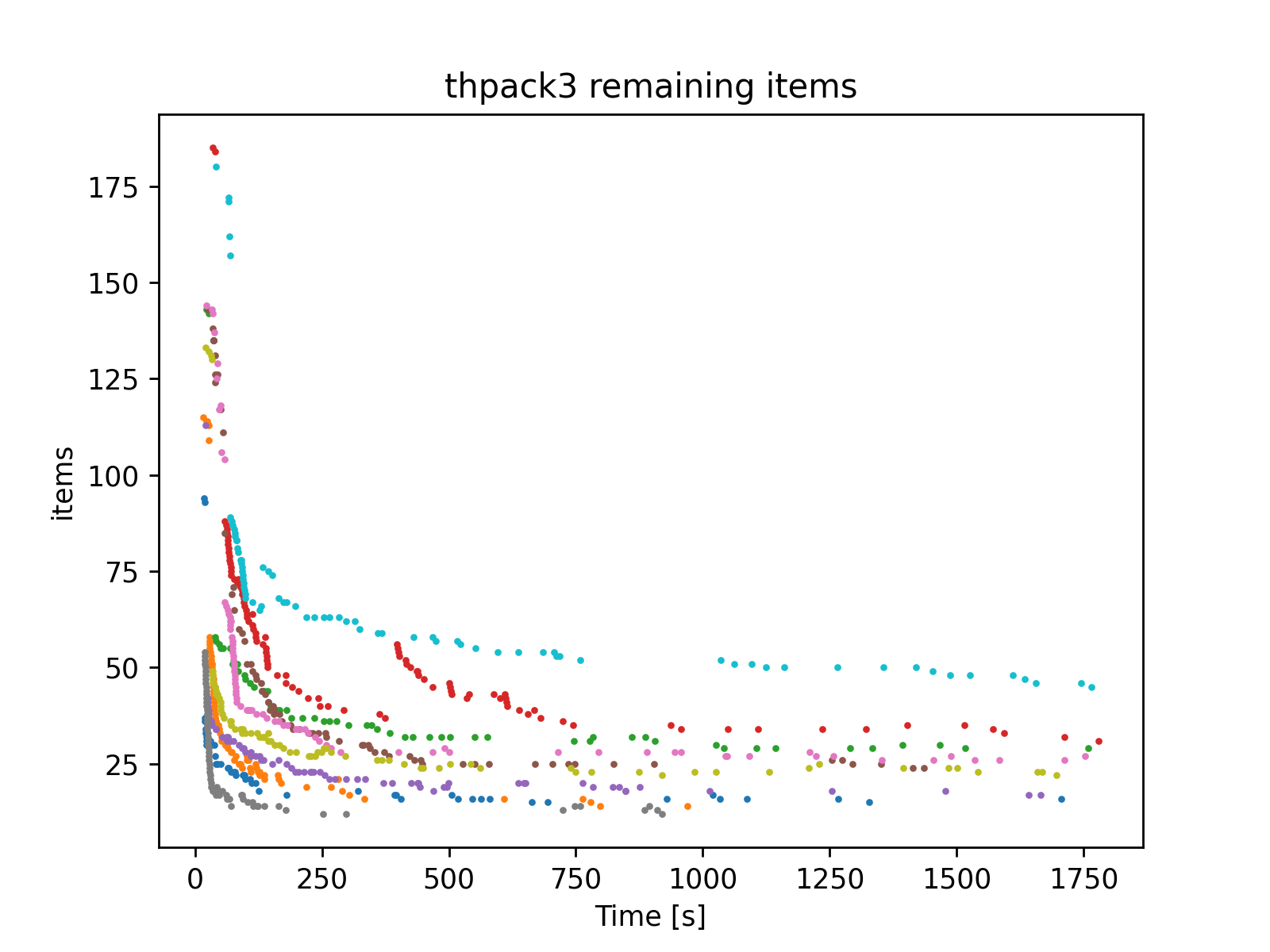}
    \caption{Remaining items \texttt{thpack3}}
    \label{fig:thpack3_boxes}    
\end{figure}

\begin{figure}
    \centering
    \includegraphics[width=0.75\linewidth]{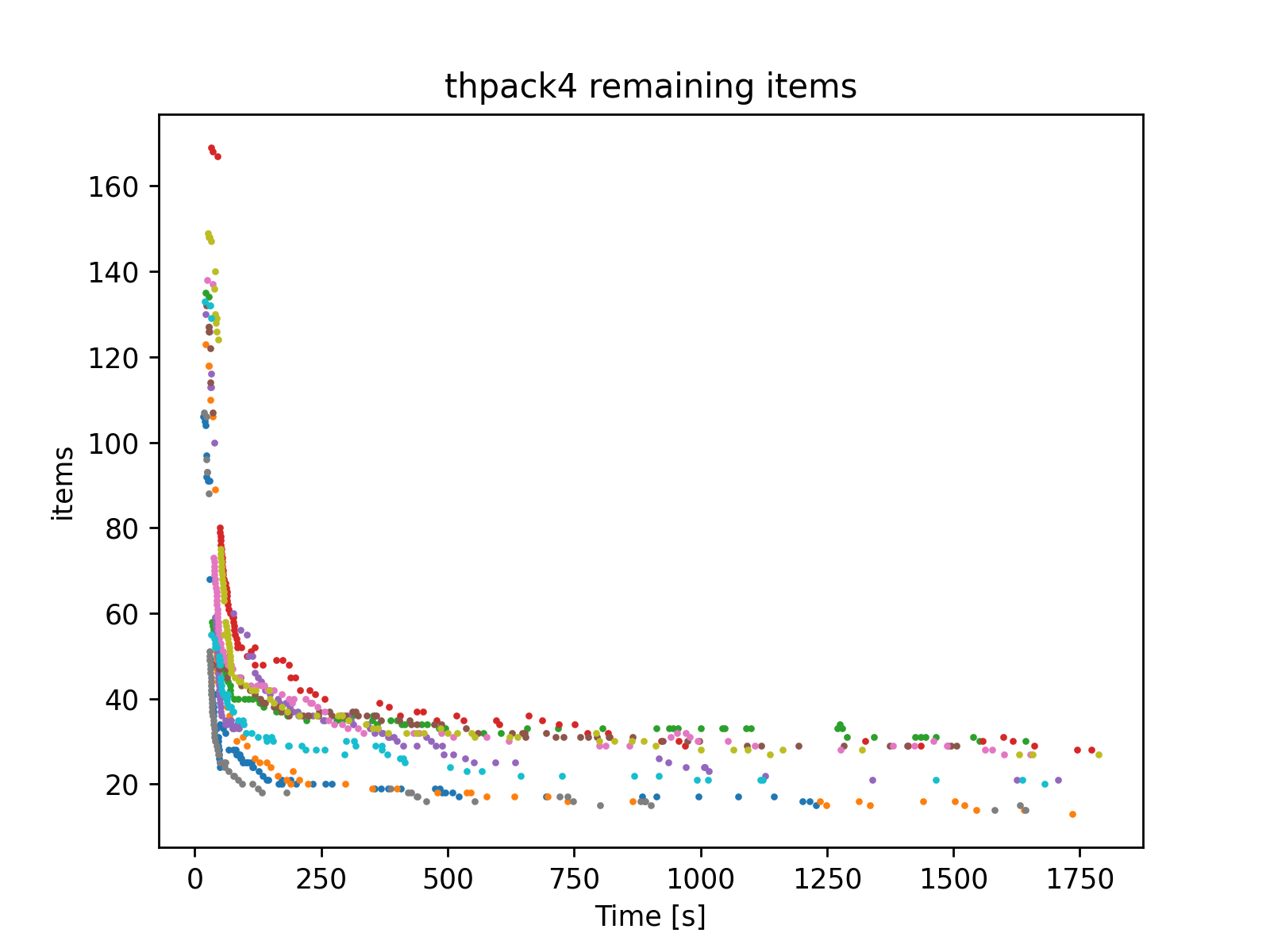}
    \caption{Remaining items \texttt{thpack4}}
    \label{fig:thpack4_boxes}    
\end{figure}

\begin{figure}[!ht]
    \centering
    \includegraphics[width=0.75\linewidth]{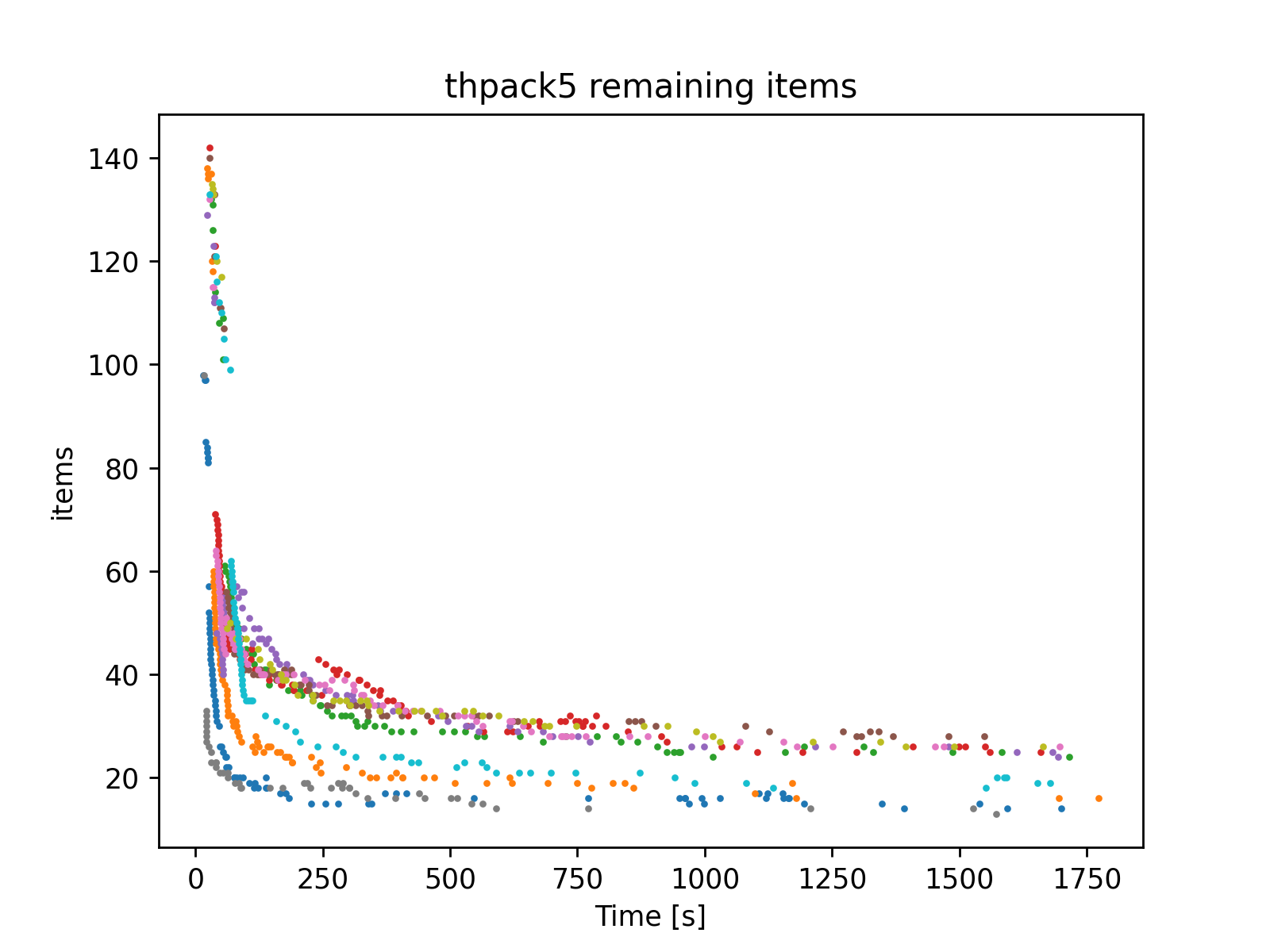}
    \caption{Remaining items \texttt{thpack5}}
    \label{fig:thpack5_boxes}    
\end{figure}

\begin{figure}[!ht]
    \centering
    \includegraphics[width=0.75\linewidth]{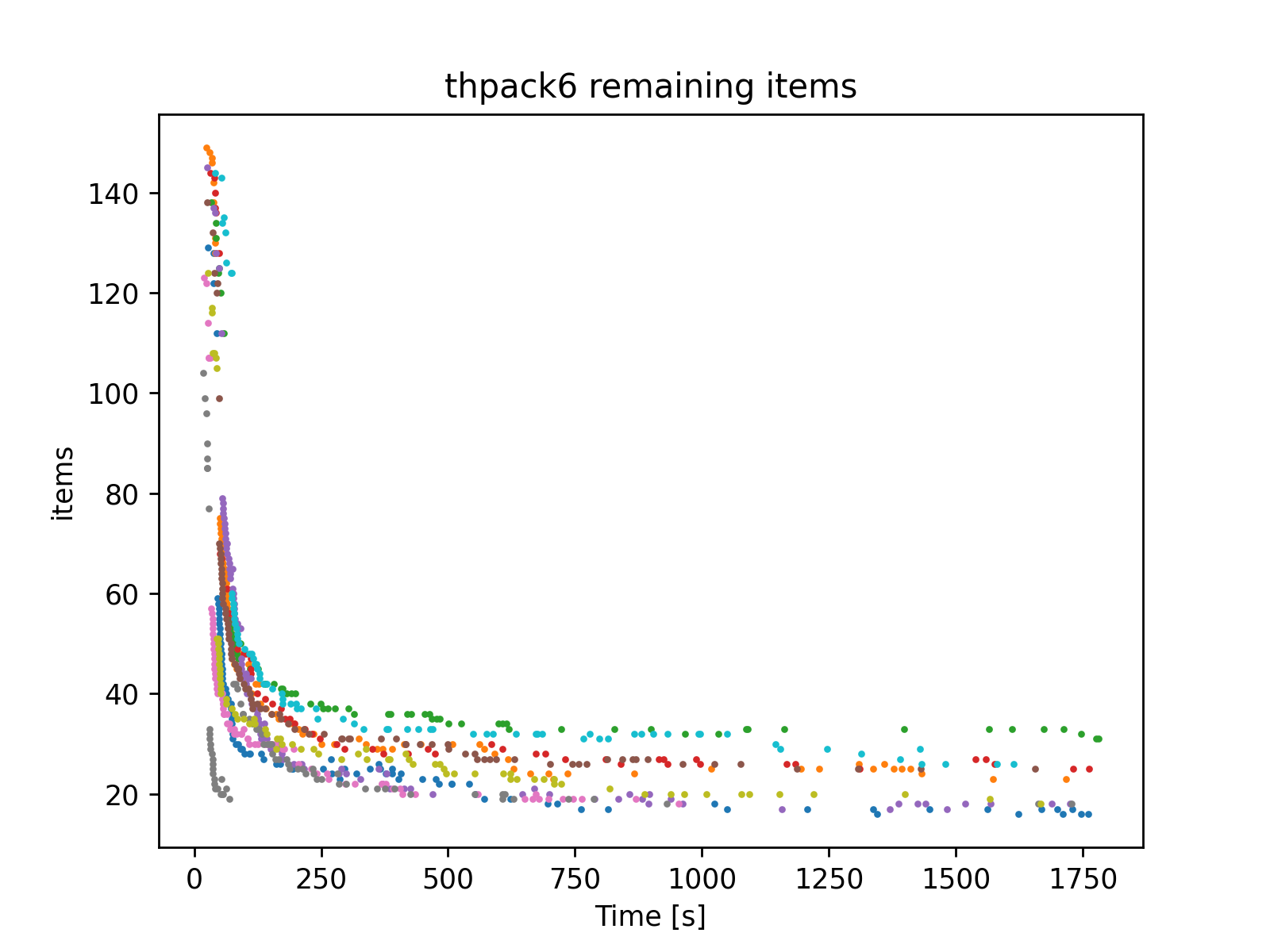}
    \caption{Remaining items \texttt{thpack6}}
    \label{fig:thpack6_boxes}    
\end{figure}

\begin{figure}[!ht]
    \centering
    \includegraphics[width=0.75\linewidth]{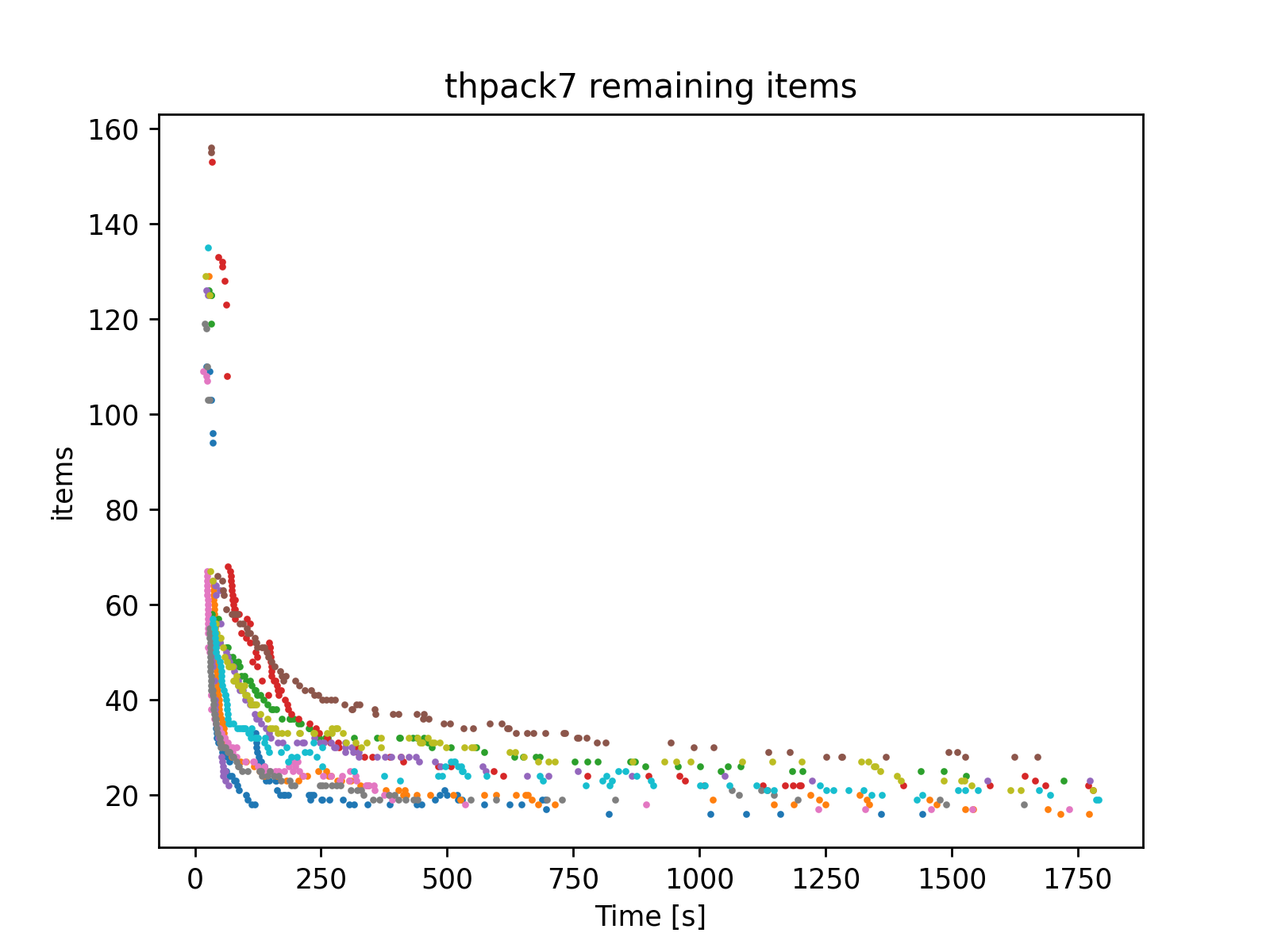}
    \caption{Remaining items \texttt{thpack7}}
    \label{fig:thpack7_boxes}    
\end{figure}

In plots \ref{fig:thpack1_boxes}-\ref{fig:thpack7_boxes}  presented is the number of remaining items in partial solutions gathered from experiments described in section \ref{sub:test-thpack1-7}.

\FloatBarrier

\subsection{Remaining payload}
\label{sec:appendix-vol}

\begin{figure}[!ht]
    \centering
    \includegraphics[width=0.75\linewidth]{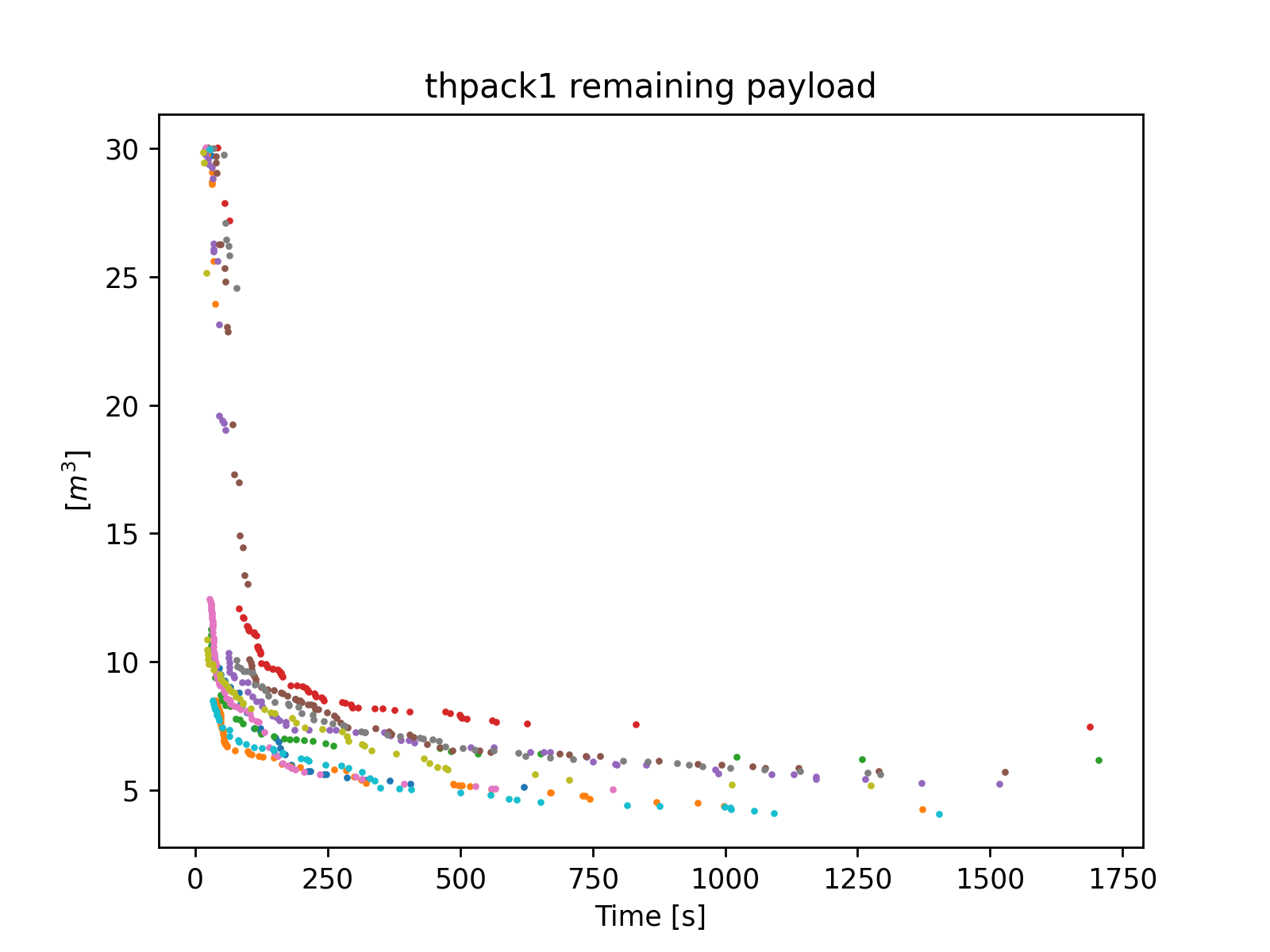}
    \caption{Remaining payload \texttt{thpack1}}
    \label{fig:thpack1_vol}    
\end{figure}

\begin{figure}[!ht]
    \centering
    \includegraphics[width=0.75\linewidth]{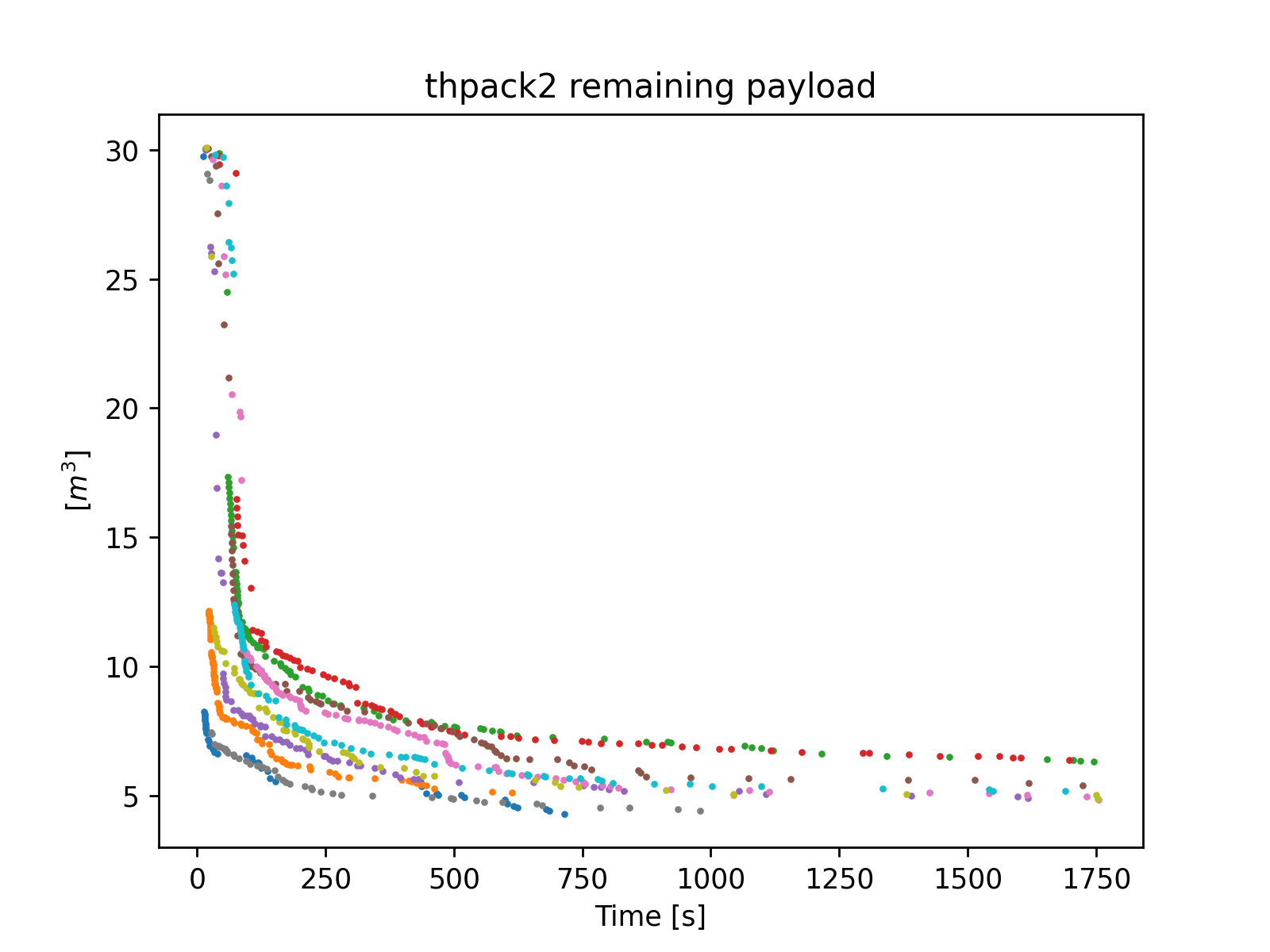}
    \caption{Remaining payload \texttt{thpack2}}
    \label{fig:thpack2_vol}    
\end{figure}

\begin{figure}[!ht]
    \centering
    \includegraphics[width=0.75\linewidth]{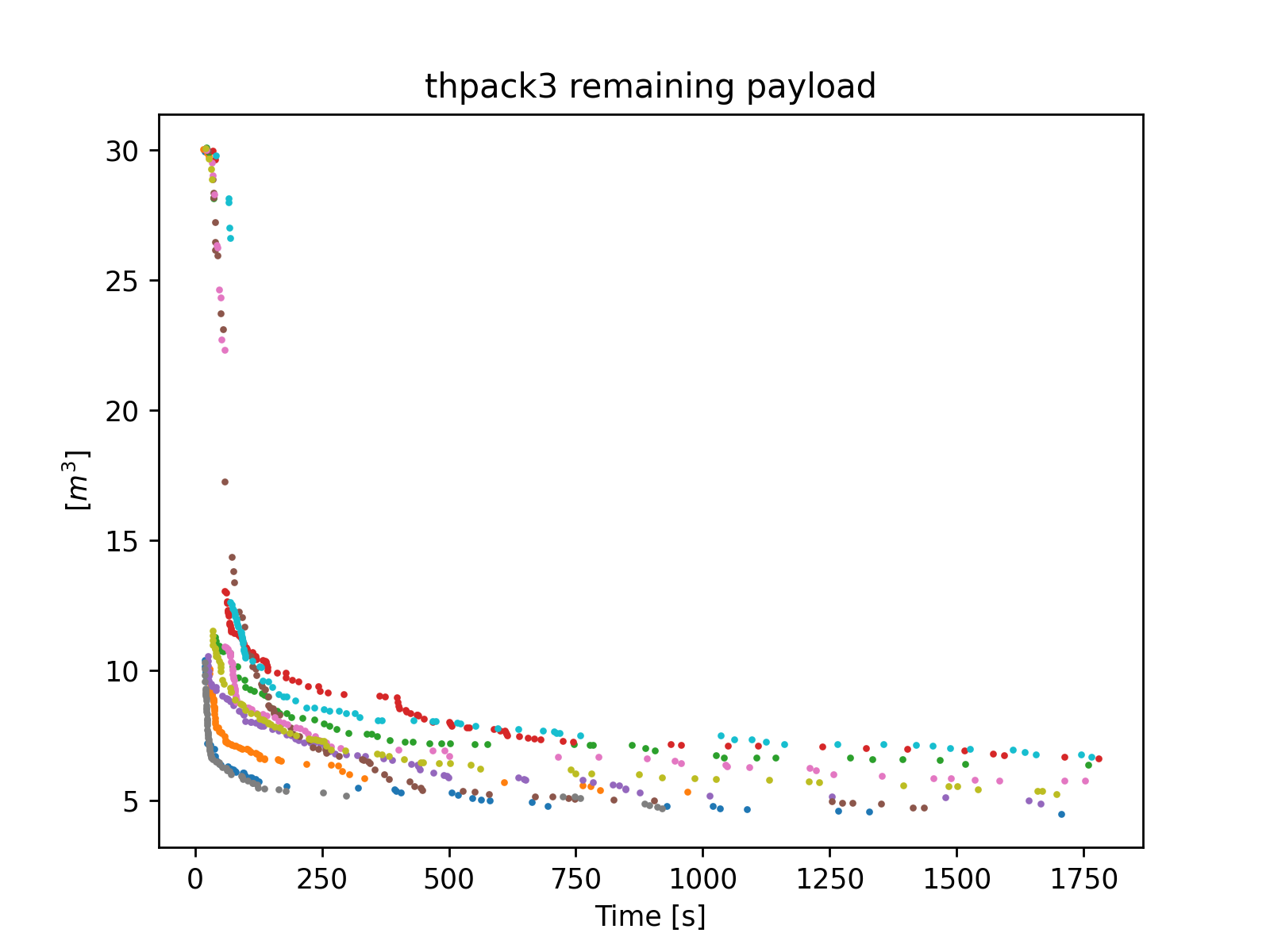}
    \caption{Remaining payload \texttt{thpack3}}
    \label{fig:thpack3_vol}    
\end{figure}

\begin{figure}[!ht]
    \centering
    \includegraphics[width=0.75\linewidth]{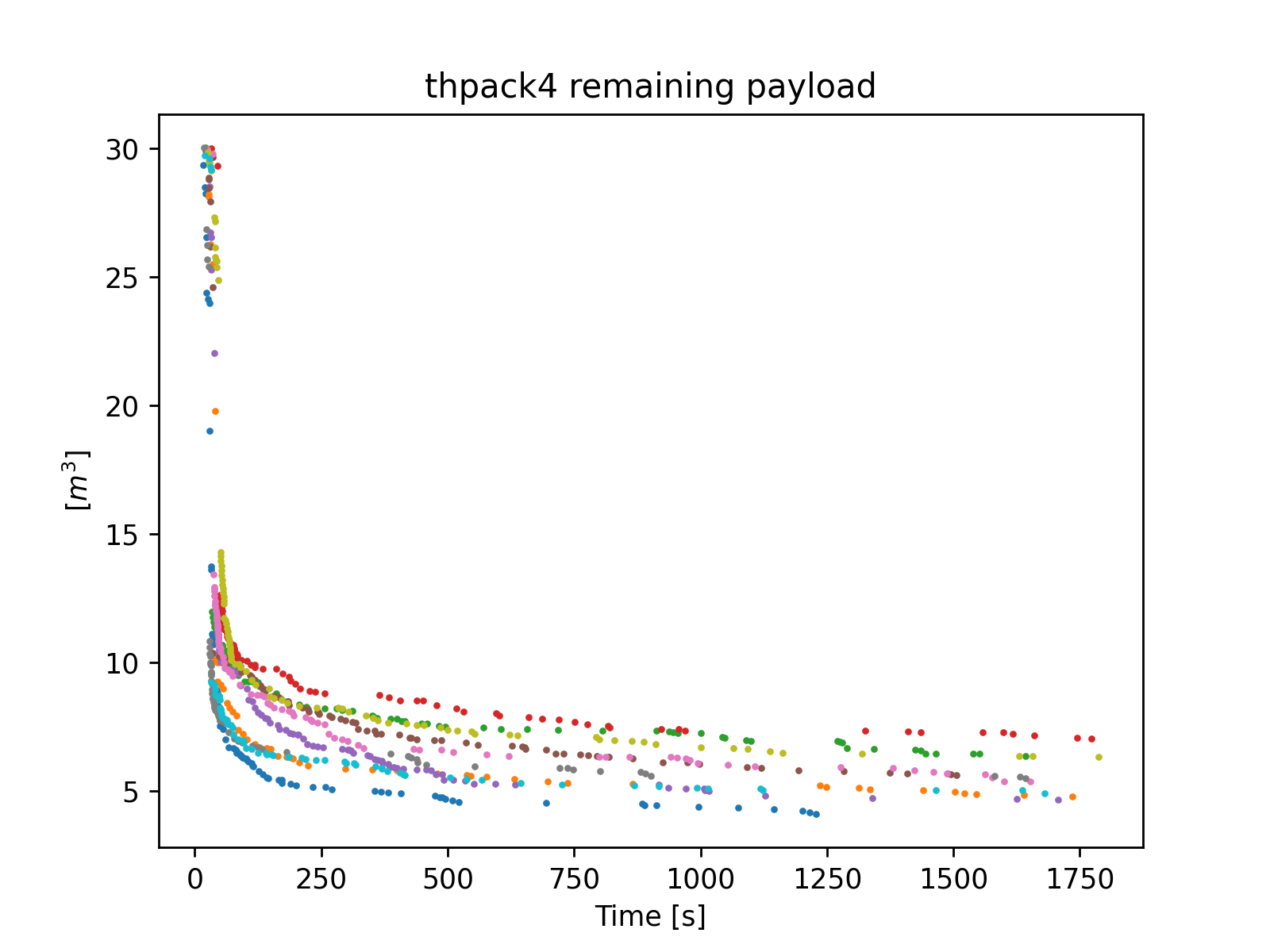}
    \caption{Remaining payload \texttt{thpack4}}
    \label{fig:thpack4_vol}    
\end{figure}

\begin{figure}[!ht]
    \centering
    \includegraphics[width=0.75\linewidth]{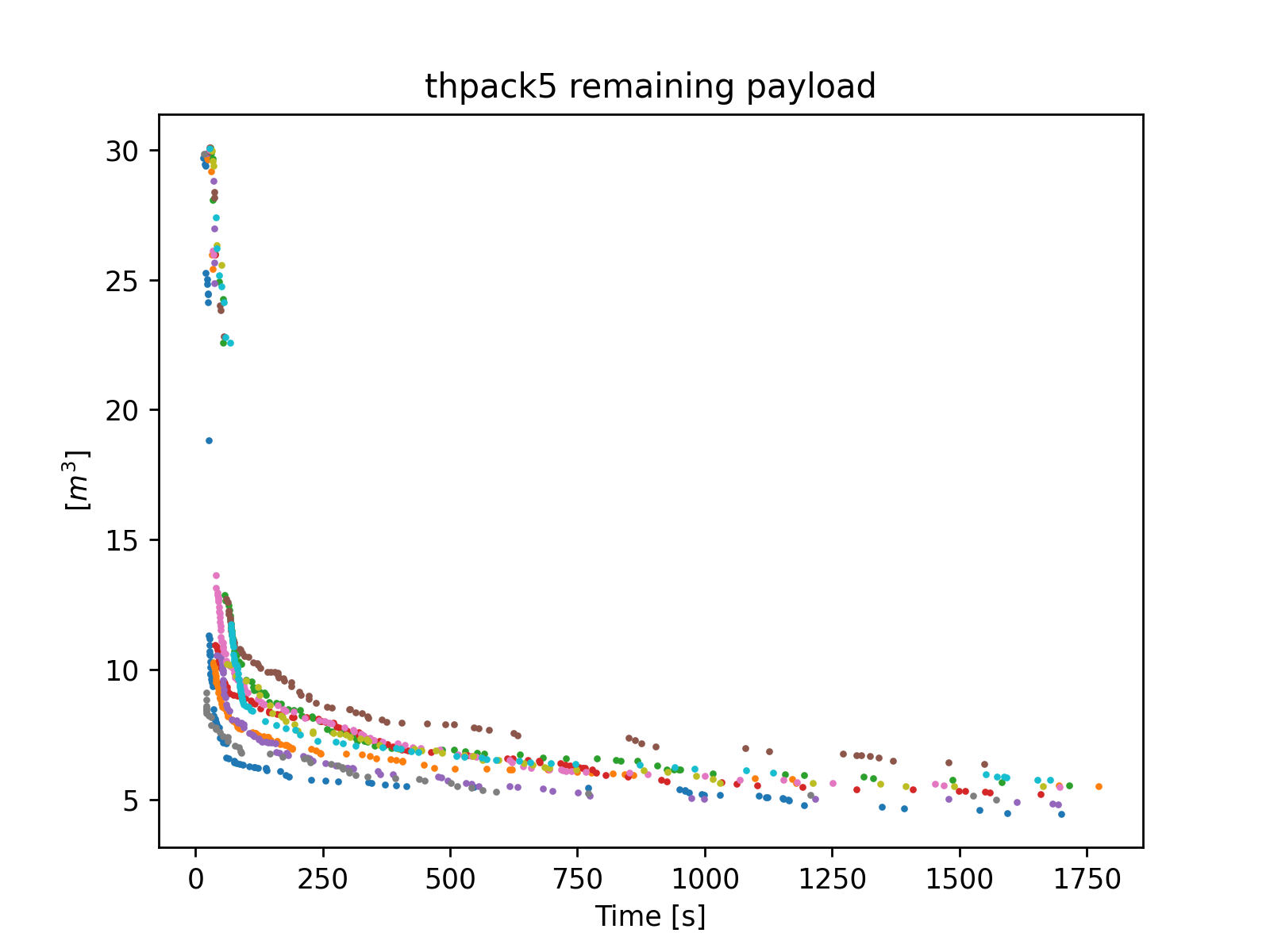}
    \caption{Remaining payload \texttt{thpack5}}
    \label{fig:thpack5_vol}    
\end{figure}

\begin{figure}[!ht]
    \centering
    \includegraphics[width=0.75\linewidth]{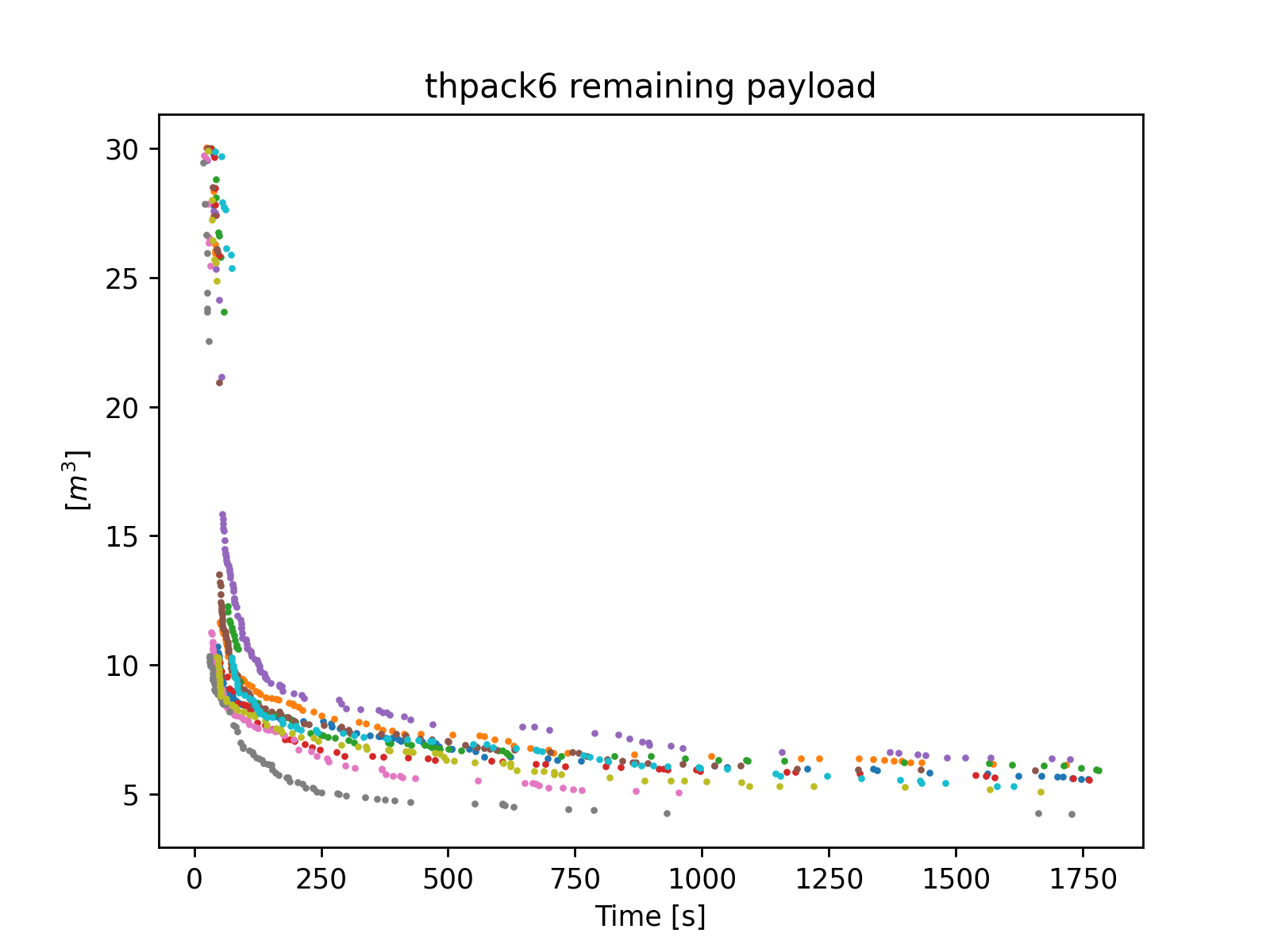}
    \caption{Remaining payload \texttt{thpack6}}
    \label{fig:thpack6_vol}    
\end{figure}

\begin{figure}[!ht]
    \centering
    \includegraphics[width=0.75\linewidth]{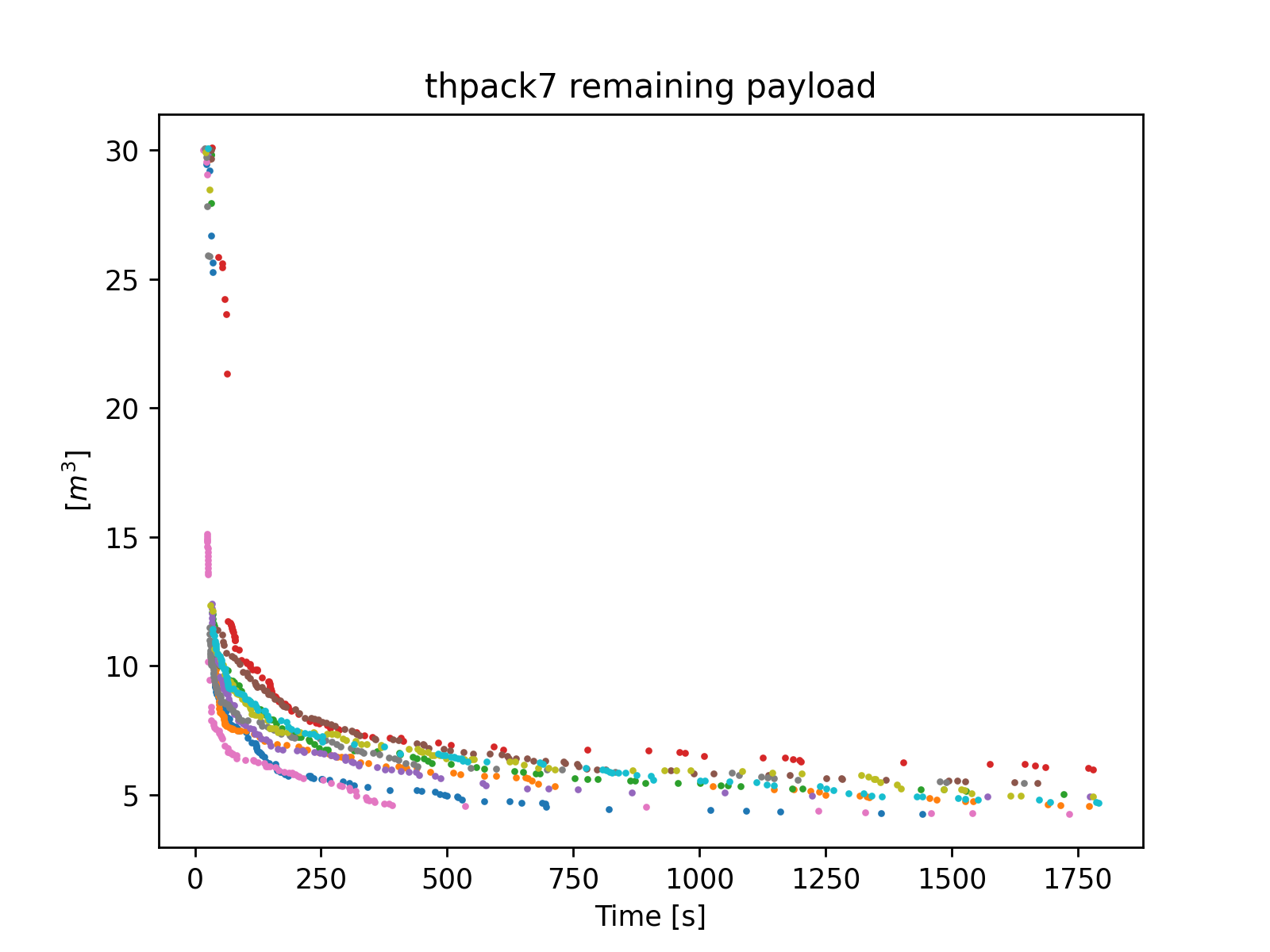}
    \caption{Remaining payload \texttt{thpack7}}
    \label{fig:thpack7_vol}    
\end{figure}

In plots \ref{fig:thpack1_vol}-\ref{fig:thpack7_vol}   presented is the volume of remaining payload in partial solutions gathered from experiments described in section \ref{sub:test-thpack1-7}.

\end{appendices}

\FloatBarrier

\bibliography{bibliography.bib}

\end{document}